%% file: Main.tex
\documentclass[twocolumn,natbib]{svjour3}

\journalname{}

\usepackage{amsmath}
\usepackage{graphicx,hyperref}
\usepackage{subcaption}
\usepackage{algpseudocode, algorithm, algorithmicx}
\usepackage{ulem, xcolor}
\usepackage{appendix}

\newcommand{\add}[0]{}
\newcommand{\delete}[1]{}
\newcommand{\addd}[0]{\color{blue}}
\newcommand{\deleted}[1]{\sout{#1}}

\begin{document}

\title{On the use of Multigrid Preconditioners for Topology Optimization}

\author{Darin Peetz \and Ahmed Elbanna}
\institute{\at Department of Civil and Environmental Engineering,\\
University of Illinois at Urbana-Champaign,\\
205 N. Mathews Ave, Urbana, IL 61801 \\
\email{peetz2@illinois.edu, elbanna2@illinois.edu}}

\maketitle

\begin{abstract}
Topology optimization for large scale problems continues to be a computational challenge.  Several works exist in the literature to address this topic, and all make use of iterative solvers to handle the linear system arising from the Finite Element Analysis (FEA).  However, the preconditioners used in these works vary, and in many cases are notably suboptimal.  A handful of works have already demonstrated the effectiveness of Geometric Multigrid (GMG) preconditioners in topology optimization. \delete{Here, we show that} {\add We provide a direct comparison of GMG preconditioners with} Algebraic Multigrid (AMG) preconditioners{\add.} \delete{offer superior robustness with only a small overhead cost.} \delete{The difference is most pronounced when the optimization develops fine-scale structural features or multiple solutions to the same linear system are needed.  We thus argue that the expanded use of AMG preconditioners in topology optimization will be essential for the optimization of more complex criteria in large-scale 3D domains.} {\add We demonstrate that AMG preconditioners offer improved robustness over GMG preconditioners as topologies evolve, albeit with a higher overhead cost. In 2D the gain from increased robustness more than offsets the overhead cost. However, in 3D the overhead becomes prohibitively large. We thus demonstrate the benefits of mixing geometric and algebraic methods to limit overhead cost while improving robustness, particularly in 3D.}
\end{abstract}

\input{Introduction.tex}
\input{Topology_Optimization.tex}
\input{Solvers.tex}
\input{Numerical_Results.tex}
\input{Discussion.tex}

\begin{acknowledgements}
The authors would like to extend their gratitude to Professor Luke Olson, UIUC Computer Science, for insightful discussions on algebraic multigrid. The authors would also like to acknowledge support from the National Science Foundation through awards 1435920 and 1753249.
\end{acknowledgements}
\begin{conflict}
The authors declare that they have no conflict of interest.
\end{conflict}
\begin{replication}
The code for this paper is available at \url{https://github.com/darinpeetz/PyOpt} and \url{https://github.com/darinpeetz/TopOpt3D}.
\end{replication}

\bibliographystyle{spbasic}
\bibliography{Mendeley}

\clearpage
\newpage
\clearpage
\appendix
{\addd
\appendixname}

\input{Appendix.tex}

\end{document}

%% file: Introduction.tex
\section{Introduction}
\label{section:Intro}
In nearly every form of continuum topology optimization, the bulk of the computational cost is incurred in either solving the linear system arising from the Finite Element Analysis (FEA) or in solving another system defined by the same linear operator (such as an adjoint equation for sensitivities). Numerous approaches have been developed in an effort to alleviate some of this cost, such as multiresolution topology optimization \citep{Kim2000Multi-resolutionParadigm,Nguyen2010AMTOP}, adaptively restricting/expanding the design space \citep{Kim2002DesignMethod}, or developing efficient, scalable methods to solve the system of equations (e.g. multigrid-preconditioned conjugate gradient \citep{Amir2014OnOptimization}). The last approach differs in that the optimization procedure itself is unaltered, and changes are isolated to only the associated finite element analysis. \delete{This paper will focus on the last approach and what improvements can remain available.} {\add This paper focuses similarly on the application of algebriac and geometric multigrid preconditioners, identifying potential limi\deleted{ta}tations and areas of improvement for each.}

There exist numerous papers in the literature exploring either 3D optimization \citep{Aage2013,Liu2014AnMatlab,Aage2015TopologyFramework,Aage2017Giga-voxelDesign} or large scale 2D optimization \citep{Amir2014OnOptimization,Jang2010DSO}, both of which generally require iterative solvers for the solution of the linear system in the finite element analysis. The efficiency of these solvers is more dependent on the choice of preconditioner than the iterative solver itself. However, even recent papers may make use of suboptimal preconditioners \citep{Benzi1999APreconditioners,Benzi2002PreconditioningSurvey} such as weighted Jacobi \citep{Mahdavi2006TopologyComputing} or incomplete Cholesky factorizations \citep{Liao2019AOptimization}. {\add These preconditioners are typically easy to set up and may be effective for sufficiently small problems, but are suboptimal in the sense that the linear solver iterations will increase as the problem size increases.} \delete{While these preconditioners are often easy to set up, their performance is less scalable than multilevel methods such as multigrid preconditioning.} Nonetheless, they still have value as smoothers within \delete{these} multilevel methods \citep{Benzi2002PreconditioningSurvey} {\add which are capable of keeping iterations constant as the problem scales up}.

A few papers have explored the use of multigrid preconditioners in topology optimization \citep{Amir2014OnOptimization,Aage2015TopologyFramework} with very promising results. However, these studies are {\add mostly} limited to only geometric multigrid (GMG) on uniform FEM grids. {\add Algebraic methods have been used in \citep{Aage2017Giga-voxelDesign}, though no discussion is provided on the relative merits of the approaches.}
To the best of our knowledge, no works exploring the {\add relative merit of AMG and GMG preconditioners} \delete{use of AMG preconditioners} in topology optimization exist in the literature. This may be at least partially attributed to the ease of implementing GMG for topology optimization. The vast majority of optimization implementations use uniform Q4 meshes in 2D or Hex8 elements in 3D. These uniform grid structures are easy to coarsen geometrically; however, the geometric approach ignores the evolution of the underlying topology.

The few examples with topology-agnostic preconditioners (such as GMG) available in the literature do not suffer a major reduction in performance as the topology evolves, but we will \delete{demonstrate} {\add describe} cases where a topology-aware preconditioner (such as algebraic multigrd, AMG) demonstrates significant improvement over the GMG approach. {\add Some strategies have also been developed to further mitigate the potential for worsened convergence with evolving topology by basing linear solver convergence on accuracy of the response function sensitivities \citep{Amir2014OnOptimization}. While this does lead to a reduction in the number of linear solver iterations, it is {\addd most applicable}\deleted{only feasible} in problems where the sensitivity calculations are cheap. In the case of compliance minimization it is very effective, but for problems {\addd like stability optimization where sensitivities depend on the solution to both an eigenvalue problem and an adjoint linear system, it is much more complicated.} \deleted{where adjoint solutions are needed it becomes entirely infeasible}. Though some of our examples would benefit from this approach, we intend to generalize our results to cases where it is not applicable and do not make use of it here.}

This paper compares the use of AMG vs. GMG preconditioners for topology optimization and {\add how their performance is influenced by the nature of the developing structure as well as algorithmic choices, specifically coarse grid size and smoother type.} \delete{what factors influence their relative performance.} We start in Section \ref{section:TO} by outlining the topology optimization framework in which the comparisons are performed. Section \ref{section:Solvers} details the methods used to solve the linear systems and generalized eigenvalue problem, including a comparison between the basic features of AMG and GMG. In Section \ref{section:Results} we present a variety of example problems to provide a numerical comparison of the performance of AMG and GMG in different scenarios. We conclude with a discussion of the findings \delete{in Section 4} and {\add in Section \ref{section:Discussion} provide} recommendations for the appropriate use of AMG and GMG preconditioners in topology optimization.

%% file: Topology_Optimization.tex
\section{Topology Optimization}
\label{section:TO}
To demonstrate the performance of the various preconditioners in topology optimization we consider two standard problems: compliance minimization and stability maximization. The first problem demonstrates the performance of the preconditioners when only a single solution to the linear system is needed and the second demonstrates the performance when multiple solutions to the same linear system are needed. In both cases we use the modified solid isotropic material with penalization (SIMP) approach \citep{Sigmund1997DesignMethod} with the linear density filter \citep{Bruns2001}.

\subsection{Compliance minimization}
The compliance minimization problem takes the following form

\begin{align}
\begin{aligned}
\min_\alpha\quad&\mathcal{F}(\boldsymbol{\alpha}) =&& \mathbf{f}^T\mathbf{u} & \\
\text{s.t.:} &&& \sum_{e=1}^{N_{el}}v_e\rho_e \leq V & \\
&&& 0 \leq \alpha_e \leq 1 & e=1,...,N_{el} \\
\text{where:} &&& \boldsymbol{\rho} = \mathbf{S}\boldsymbol{\alpha} & \\
&&& \mathbf{K}(\mathbf{\rho})\mathbf{u} = \mathbf{f} &
\end{aligned}
\end{align}
where $\mathbf{\alpha}$ represents the design variables, $\mathbf{S}$ is the filtering matrix, $\mathbf{\rho}$ represents the filtered densities, $\mathbf{K}$ is the stiffness matrix, $\mathbf{u}$ is the vector of displacements, $\mathbf{f}$ is the vector of external forces, $v_e$ is the volume of element $e$, and $V$ is the total allowable volume of the structure. To prevent the stiffness matrix from becoming singular, the element stiffnesses used to construct the local stiffness matrices are calculated using the modified SIMP rule:
\begin{equation}
    E(\rho) = E_{min} + (E_{\max} - E_{min})\rho^p
\end{equation}
where the penalty, $p$ is gradually increased from 1 to 4 and $E_{min}$ is set as $1e-10E_{\max}$.

Using the adjoint method, the sensitivities of the objective function are calculated as
\begin{align}
\begin{aligned}
    \frac{\partial \mathcal{F}}{\partial \rho_e} &= -\frac{\partial E_e}{\partial \rho_e}\mathbf{u}^T\frac{\partial \mathbf{K}}{\partial E_e}\mathbf{u} \\
    \frac{\partial \mathcal{F}}{\partial \alpha} &= \mathbf{S}^T\frac{\partial \mathcal{F}}{\partial \rho}
\end{aligned}
\end{align}
For the sake of generality we will use the method of moving asymptotes (MMA) \citep{Svanberg1987} for design updates.

\subsection{Stability maximization}
The problem of optimizing for structural stability takes the following form

\begin{align}
\begin{aligned}
\min_\alpha\quad&\mathcal{F}(\boldsymbol{\alpha}) =&& \frac{1}{P_{critical}} = \lambda_{\max} & \\
\text{s.t.:} &&& \sum_{e=1}^{N_{el}}v_e\rho_e \leq V & \\
&&& 0 \leq \alpha_e \leq 1 & e=1,...,N_{el} \\
\text{where:} &&& \boldsymbol{\rho} = \mathbf{S}\boldsymbol{\alpha} & \\
&&& \mathbf{K}(\boldsymbol{\rho})\mathbf{u} = \mathbf{f} & \\
&&& \mathbf{K}_\sigma\mathbf{\Phi} = \lambda\mathbf{K}\mathbf{\Phi}&
\end{aligned}
\end{align}
where $\mathbf{K}_\sigma$ is the stress stiffness matrix, $\mathbf{\Phi}$ is the eigenvector of the generalized system, and $\lambda$ is the corresponding eigenvalue. Because $\mathbf{K}_\sigma$ is potentially indefinite, it is most natural to write the generalized eigenvalue equation in this form so that the matrix on the right-hand-side is positive definite, which is assumed for many eigenvalue solvers. To prevent critical buckling modes from appearing in non-structural regions of the domain we use another modified version of the SIMP formula to interpolate the values of stiffness for the stress stiffness matrix  \citep{Bendse2003,Gao2015TopologyConstraints,Thomsen2018BucklingAnalysis}

\begin{equation}
    E_\sigma(\rho) =
    \begin{cases}
        E_{\max}\cdot\rho^p, & \text{if } \rho >= 0.1 \\
        0, & \text{if } \rho < 0.1
    \end{cases}
\end{equation}

Again using the adjoint method, we can derive the sensitivities of the stability problem as
\begin{align}
\begin{aligned}
    \frac{\partial \mathcal{F}}{\partial \rho_e} &= \mathbf{\Phi}\left(\frac{\partial E_{\sigma,e}}{\partial \rho_e}\frac{\partial \mathbf{K}_{\sigma,e}}{\partial E_e} - \lambda_{\max}\frac{\partial E_e}{\partial \rho_e}\frac{\partial \mathbf{K}}{\partial E_e}\right)\mathbf{\Phi} \\
    & \quad\quad\quad + \mathbf{v}^T \frac{\partial E_{\sigma,e}}{\partial \rho_e}\frac{\partial \mathbf{K}}{\partial E_e} \mathbf{u}\\
    \frac{\partial \mathcal{F}}{\partial \boldsymbol{\alpha}} &= \mathbf{S}^T\frac{\partial \mathcal{F}}{\partial \boldsymbol{\rho}}
    \label{equ:Stability_Sensitivity}
\end{aligned}
\end{align}
where $\mathbf{v}$ is the solution to the adjoint equation defined by
\begin{equation}
    \mathbf{K}\mathbf{v} = \mathbf{\Phi}^T\frac{\partial \mathbf{K}_\sigma}{\partial \mathbf{u}}\mathbf{\Phi}
\end{equation}
Here we see that the adjoint equation requires another solution to a linear system defined by $\mathbf{K}$ for each eigenvalue calculated, in addition to the multiple solutions required inside any solver for the generalized eigenvalue problem.

While the formulation shown here includes only the maximal eigenvalue in the optimization, in practice it is necessary to aggregate a subset of the largest eigenvalues \delete{\citep{Ferrari2019RevisitingConstraints}}. {\add Perhaps the most important motivation for aggregation is the fact that eigenvalues with algebraic multiplicity greater than one have undefined sensitivities \citep{Seyranian1994}. However, if all of the multiple eigenmodes are aggregated, the sensitivity of the aggregate function becomes well-defined even in the case of multiplicity \citep{Ferrari2019RevisitingConstraints,Gravesen2011OnEigenvalues,Torii2017StructuralApproximation}. A secondary motivation is that aggregation also provides greater continuity if eigenmodes switch order (or simply switch between single and multiple eigenvalues) as part of the optimization process.}

\delete{For the purposes of this paper we optimize the 1-norm of the set of $n_{eig}$ largest eigenvalues such that any two consecutive eigenvalues in this set are separated by no more than 1\% of each other, possibly changing the number of eigenvalues considered at every iteration. With these considerations in mind, the new objective and its sensitivities become:}
{\add For this paper we use the $p$-norm strategy described in \citep{Ferrari2019RevisitingConstraints}. Specifically we take the objective of the stability function to be the $p=8$-norm of the first 6 eigenmodes. This changes the stability function and sensitivity to be:

\begin{equation}
    \min_\alpha\quad\mathcal{F}(\boldsymbol{\alpha}) = \left(\sum_{i=1}^{n=6} \lambda_{i}^8\right)^{1/8}
\end{equation}
\begin{align}
\begin{aligned}
    \frac{\partial \mathcal{F}}{\partial \rho_e} &= \mathcal{F}^{(1-8)}\sum_{i=1}^n \lambda_{i}^{(8-1)}\frac{\partial \lambda_i}{\partial \rho_e}
\end{aligned}
\end{align}
where $\Phi_i$, $\lambda_i$, and $v_i$ are calculated independently for each of the $n$ modes clustered near the end of the spectrum. Note that $\frac{\partial \lambda_i}{\partial \rho_e}$ shares its definition with $\frac{\partial \mathcal{F}}{\partial \rho_e}$ in Equation \ref{equ:Stability_Sensitivity}.}

%% file: Solvers.tex
\section{Solvers}
\label{section:Solvers}
\subsection{Linear Solvers}
\label{section:Linear_Solvers}
First we discuss the methods available to solve a single linear system, either for calculating displacements or solving the adjoint equation. The system of equations takes the form $\mathbf{K}\mathbf{U} = \mathbf{F}$, though $\mathbf{U}$ and $\mathbf{F}$ may represent something other than the displacements and external forces in the case of the adjoint equation. $\mathbf{K}$ is an arbitrarily large, sparse, symmetric, and positive definite matrix with an approximately constant number of nonzeros per row. In this case, the sparse Cholesky factorization operates in $\mathcal{O}(n^{3/2})$ time ($n$ being the size of the stiffness matrix) \citep{Davis2006DirectSystems}. For small 2D problems this scaling is generally satisfactory for solving the linear system, and offers the additional bonus that subsequent solutions to the same linear system can make use of the same factorization and operate in closer to $\mathcal{O}(n)$ time.

The advantages of the Cholesky factorization, and direct solvers in general, begin to fade for large problems in 2D and even for smaller problems in 3D. The factorization requires substantially more memory than the matrix itself and as $n$ grows, the difference in $\mathcal{O}(n^{3/2})$ and $\mathcal{O}(n)$ time becomes significant. In these cases, iterative solvers \citep{Saad2003} are more attractive as their cost is dominated by $\mathcal{O}(n)$ matrix-vector operations, and storage requirements beyond the matrix itself are limited to a small set of vectors of length $n$. The Conjugate Gradient method \citep{Hestenes1952} is {\add often} preferred for sparse, symmetric, positive-definite (SPD) matrices because it most effectively makes use of these properties{\add.} \delete{, as opposed to more general methods like Generalized Minimum Residual (GMRES) \citep{Saad1986}.} {\add Other, more general methods like Generalized Minimum Residual (GMRES) \citep{Saad1986} are less efficient in theory but often perform as well or better in practice, depending on the nature of the linear system.} \delete{While the CG method is guaranteed to converge in $n$ iterations, in practice a certain amount of error is allowed in the solution and iterations are stopped when a convergence tolerance is reached. Thus,} {\add Regardless of the particular method used}, the effectiveness of the method for solving a linear system is governed (roughly) by the convergence rate
\begin{equation}
    \frac{\sqrt{\kappa}-1}{\sqrt{\kappa}+1}
\end{equation}
where 
\begin{equation}
    \kappa = \text{cond}(\mathbf{K}) = \frac{\lambda_{\max}(\mathbf{K})}{\lambda_{\min}(\mathbf{K})}
\end{equation}

In the case of topology optimization $\kappa \gg 1$, especially after the optimization begins to produce a structure, and void regions of the design domain assume a stiffness several orders of magnitude smaller than the solid regions. For a fully converged solid-void structure, the lowest energy modes are almost always confined to the void regions while the high energy modes exist primarily in the solid regions, meaning that $\kappa$ \delete{scales roughly with $E_{min}$} {\add worsens as $E_{min}$ decreases}. {\add Approaches such as design space optimization \citep{Kim2002DesignMethod} that remove these regions from the linear system may improve conditioning, but run the risk of isolating structural elements, thereby creating a nullspace. Even if this phenomenon is avoided, the varying topology will still produce very poorly-conditioned linear systems.}

\delete{The ill-conditioning of the system is not unique to topology optimization, but is instead very common when using iterative solvers for any linear system.} The poor conditioning is overcome through the use of a quality preconditioner, $\mathbf{M}$. The preconditioner, which may itself be another matrix or simply an operator, is chosen so that $\kappa(\mathbf{M}^{-1}\mathbf{K}) \ll \kappa(\mathbf{K})$. Most of the literature on topology optimization uses direct methods to solve the linear system, where $\mathbf{M} = \mathbf{K}$ explicitly. There are a handful of works that explore the use of other preconditioners, particularly for optimizations in 3D \citep{Aage2013,Aage2015TopologyFramework} or with a large number of degrees of freedom in 2D \citep{Mahdavi2006TopologyComputing}. Some of these works have shown great promise using multigrid preconditioners \citep{Amir2014OnOptimization,Kennedy2015Large-scaleManufacturing}, though they focus exclusively on geometric versions. {\add One notable exception is \citep{Aage2017Giga-voxelDesign}, which makes use of a combination of geometric and algebraic methods, though it offers no insight on the relative merits of each, as we do here.}

\subsection{Multigrid Methods}
\label{section:MG}
Multigrid methods are a class of multilevel methods that replicate the discretization of the original partial differential equation (PDE) on increasingly coarser grids to improve performance. They are based on the premise that while smoothers such as Gauss-Seidel or Jacobi may be ineffective at reducing all of the error in a solution approximation, they are very effective at removing error with large eigenvalues (high-energy errors) \citep{Saad2003,Briggs2000}. {\add These error modes typically appear as highly oscillatory modes on the given discretization.} When the remaining ``smooth" errors are projected onto a coarser grid, {\add a fraction of them will be represented as oscillatory modes in the new system.} \delete{they in turn correspond to high-energy errors in the new system matrix.} \delete{When the original error is recursively smoothed and projected onto coarser and coarser grids,} {\add Through recursive projection most of the modes become oscillatory and can in turn be removed with smoothing.} Eventually the problem becomes small enough to make use of a direct solver {\add for any remaining error modes}. In this way, the operations that take place on the full original grid are limited to matrix-vector operations and the more expensive operations are performed on smaller grids where the cost becomes negligible.

{\add It is also possible to reduce the amount of coarsening needed by using inexact or iterative solvers on the coarsest grid. In such a scenario some error associated with the smoothest modes of the operator may be omitted in the correction, but the small eigenvalues associated with these modes mean that correcting them offers little improvement to the residual on the fine grid anyway. Though not explored here, this approach has been used successfully in other works \citep{Aage2015TopologyFramework}.}

\begin{figure}
	\centering
	\includegraphics[width=0.45\linewidth]{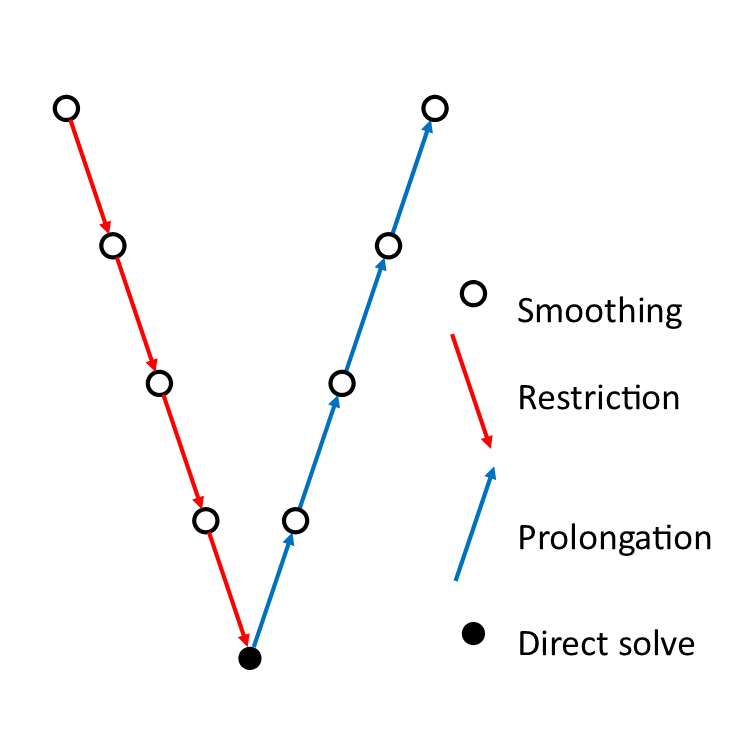}
	\caption{Multigrid V-Cycle}
	\label{fig:Multigrid}
\end{figure}

A typical V-Cycle multigrid algorithm is presented in Algorithm \ref{alg:MG} and \delete{visually} {\add visualized} in Figure \ref{fig:Multigrid}. The procedure is the same regardless of whether an algebraic or geometric method is used to assemble the grid hierarchy. Geometric methods, which have already been used for topology optimization in \citep{Aage2013,Amir2014OnOptimization,Aage2015TopologyFramework} {\add among others}, construct the hierarchy directly from the physical grid. \delete{The coarse grid in this case}{\add The restriction/prolongation operators used to transfer vectors between the levels in the hierarchy are defined directly by the shape functions of a coarser discretization of the original mesh. In the case of uniform material density, this means the coarse grid} is essentially a rediscretization of the original PDE with half of the number of elements in each direction. \delete{The} {\add This definition of the} prolongation operator, $\mathbf{P}$ \delete{is then constructed from the finite element interpolation on the coarse mesh so} {\add ensures} that the solution on the coarse grid is reconstructed exactly on the fine grid. As the system matrix $\mathbf{A}^0$ ($\mathbf{K}$ in the case of topology optimization) is symmetric, we use Galerkin projections and set the restriction matrix $\mathbf{R} = \mathbf{P}^T$. Prior to solving the linear system, a series of linear operators is constructed, one for each level of the multigrid. The operator on the fine grid is given from the original linear system of equations for displacements, and each subsequent operator is defined by the Galerkin projection as
\begin{equation}
    \mathbf{A}^{i+1} = \mathbf{P}^T\mathbf{A}^i\mathbf{P}
    \label{equ:Galerkin_Projection}
\end{equation}
If a direct solver is used on the coarsest level, the coarse operator should also be factorized at this time.

\begin{algorithm}[hbt]
	\caption{V-cycle multigrid with weighted Jacobi smoothing for the system $\mathbf{A}\mathbf{x}=\mathbf{b}$}
	\label{alg:MG}
	\begin{algorithmic}
		\Procedure{MG}{$\mathbf{A}, \mathbf{P}, \mathbf{b}, k, n_{levels}, n_{pre}, n_{post}, w$}
        \State $\mathbf{x} = \mathbf{0}$
        \For{$i=1..n_{pre}$}
        \State $\mathbf{r} = \mathbf{b} - \mathbf{A}^k\mathbf{x}$
        \State $\mathbf{x} += w*(\text{diag}(\mathbf{A}^k)^{-1}\mathbf{r})$
        \EndFor
        \State $\hat{\mathbf{b}} = (\mathbf{P}^k)^T(\mathbf{b}-\mathbf{A}^k\mathbf{x})$
        \If{$k<n_{levels}-1$}
        \State $\mathbf{x} += \mathbf{P}^k*\text{MG}\left(\mathbf{A},\mathbf{P},\hat{\mathbf{b}}, k+1, n_{levels}, n_{pre}, n_{post}, w\right)$
        \Else
        \State $\mathbf{x} += \mathbf{P}^k*(\mathbf{A}^k)^{-1}\hat{\mathbf{b}}$
        \EndIf
        \For{$i=1..n_{post}$}
        \State $\mathbf{r} = \mathbf{b} - \mathbf{A}^k\mathbf{x}$
        \State $\mathbf{x} += w*(\text{diag}(\mathbf{A}^k)^{-1}\mathbf{r})$
        \EndFor
        \State \Return{$\mathbf{x}$}
		\EndProcedure
	\end{algorithmic}
\end{algorithm}

Algebraic methods differ from geometric methods in that the levels in the multigrid are constructed not based on the mesh, but instead are constructed directly from the linear operator itself. This gives the method flexibility to be applied to problems where a regular mesh may not be available, or grid restriction does not accurately capture smooth modes of the operator (as in the case of anisotropic diffusion \citep{Saad2003}). There are various methods available to perform this hierarchy assembly, going back to the classical Ruge-St\"uben \citep{Ruge19874.Multigrid}. All methods follow a general procedure of identifying which nodes (or degrees of freedom) are strongly connected to each other and lumping them into ``supernodes" as part of the restriction operation. The definition of ``strongly connected" is somewhat heuristic, but in the case of topology optimization it ensures that nodes attached to structural elements are not lumped with nodes in the void regions. For this work we will focus on the smoothed aggregation method \citep{Vanek1996AlgebraicProblems} due to its superior performance and wide availability in software packages.

The assembly of the operator hierarchy in AMG works similarly to GMG, but with a few additional steps. Starting with the fine grid operator $\mathbf{A}^0$, the off-diagonal elements of the matrix are compared to diagonal elements to determine which degrees of freedom are strongly connected. A restriction operator, $\mathbf{R}$ is then constructed based on the type of AMG used, which projects smooth error from the fine system onto a coarsened system where multiple strongly connected degrees of freedom are represented with only a few reduced degrees of freedom. {\add As in the case of the geometric approach, we take advantage of the symmetry of the system to define coarse grid operators using Equation \ref{equ:Galerkin_Projection}. Once the entire hierarchy is constructed any additional preparations (such as factorizing the coarsest operator) are performed.} \delete{Again taking advantage of the symmetry of the system and using the Galerkin projection to define coarse grid operators using Equation \ref{equ:Galerkin_Projection} and factorize the coarsest operator.}

\subsection{Eigenvalue solvers}
\label{section:Eigensolvers}
The performance of eigenvalue solvers, particularly for the generalized eigenvalue problem, is more complicated than that of linear solvers. While all practical eigensolvers must be iterative methods \citep{Abel}, generalized eigensolvers internally require a solution to a linear system at each iteration, which may be solved using any of the previously described methods for solving linear systems of equations. Some methods, such as Arnoldi or Lanczos, \citep{ARNOLDI} make use of a factorization of one of the matrices ($\mathbf{K}$ in this case as it is positive definite) to get a high-precision solution each time. Other methods, so-called preconditioned eigensolvers, make use of an iterative method to solve (or approximate the solution to) a linear system. As in the case of simple linear systems, matrix factorizations are really only feasible for smaller problems and the preconditioned eigensolvers are necessary for large-scale problems. It should be noted that in the eigenvalue problem only one factorization is needed, but the system must be solved with many right-hand-sides. This suggests that the tradeoff in efficiency between factorized and preconditioned methods occurs at a larger system size in the eigenvalue problem than that of the simple linear system.

Of the preconditioned eigensolvers, three are widely used and demonstrate good performance in the problem type we are examining. Two are variations of Davidson's method \citep{Morgan1986GeneralizationsMatrices}: generalized Davidson \citep{Davidson1975TheMatrices} and Jacobi-Davidson \citep{G.Sleijpen1996,Sleijpen1996}, and the other is the locally optimal block preconditioned conjugate gradient (LOBPCG) \citep{Knyazev2001}. While all three methods use a preconditioner to solve a linear system, they differ in how that system is defined and in how accurately it needs to be solved. The Jacobi-Davidson method, with its internal correction equation, is more suited to target interior eigenvalues \citep{Sleijpen1996}, whereas generalized Davidson and LOBPCG are more effective for problems where exterior eigenvalues are needed (as in the case of stability optimization). In this paper we use the generalized Davidson's method {\add (described in more detail in Algorithm \ref{alg:GD})} because the larger search space it uses gives it better performance in the context of stability optimization than LOBPCG. {\add For the examples shown, we use a minimum and maximum search space size ($j_{min}, j_{max}$) of 10 and 25, respectively.}

\begin{algorithm}[hbt]
    \add
	\caption{Generalized Davidson Eigensolver}
	\label{alg:GD}
	\begin{algorithmic}[1]
		\Procedure{GD}{$\mathbf{A}, \mathbf{B}, M$}
		\State Initialize search space $\mathbf{V} = \{\mathbf{v}_1, \mathbf{v}_2, ..., \mathbf{v}_{j_{min}}\}$
		\State Orthonormalize $\mathbf{V}$ using Gram-Schmidt
		\State $it = 0$
		\While{$it < it_{max}$}
		\State --- Projected eigenproblem ---
		\State Compute eigenmodes $(\mathbf{S}, \mathbf{Q})$, on span of $V$ using Rayleigh-Ritz
		\State Pick closest mode to target eigenvalue, $\theta, \mathbf{q}$
		\State $\mathbf{r} = \mathbf{A}\mathbf{q} - \theta \mathbf{B}\mathbf{q}$
		\State --- Reset search space if necessary ---
		\If{$j = j_{max}$}
		\State $\mathbf{V} = \mathbf{V} * \{\mathbf{q}_1, \mathbf{q}_2, ..., \mathbf{q}_{j_{min}}\}$
		\EndIf
		\State --- Expand search space ---
		\State $\mathbf{z} = M(\mathbf{r})$ (apply preconditioner to residual)
		\State Orthonormalize $\mathbf{z}$ against $\mathbf{V}$
		\State $\mathbf{V} = \{\mathbf{V}, \mathbf{z}\}$
		\State $j += 1$
		\EndWhile
		\EndProcedure
	\end{algorithmic}
\end{algorithm}

%% file: Numerical_Results.tex
\section{Numerical Results}
\label{section:Results}
To demonstrate the relative performance of AMG and GMG we will look at \delete{four} {\add a variety of} different topology optimization cases {\add in 2 and 3 dimensions}. {\add The majority of the problems are run in parallel where we make use of the PETSc library \citep{petsc-efficient,petsc-user-ref,petsc-web-page} to construct the AMG hierarchy and to provide the interface for applying both the AMG and GMG preconditioners as part of a linear solver. For the grid problem (which is unique in that it is not run in parallel) we use the PyAMG library \citep{OlSc2018} to provide the same functionality. To construct the GMG hierarchy we use a Galerkin projection with restriction/prolongation operator coefficients defined by shape function values on the coarse grids, as used in \citep{Amir2014OnOptimization}.} \delete{For every case the AMG and GMG preconditioners are set up with enough levels that the coarsest grid has no more than 80 nodes. At each level of the preconditioner a single pre- and post-smoothing pass is applied using weighted Jacobi, and the coarsest level is solved with an LU decomposition. Within the AMG preconditioner we also compare the performance of a blocked (one block per node) and unblocked weighted Jacobi smoother. Our results indicated no difference in performance between blocked and unblocked versions for the GMG preconditioner, so only the results for the block smoother are shown here. For brevity, the three schemes will be referred to as GMG, AMG, and block AMG for the remainder of this manuscript.}

For the AMG preconditioner, strength of connection is calculated in a block fashion (one value per node) using the symmetric strength of connection formulation with {\add$\beta = 0.003$
\begin{equation}
    S = \mathbf{K}(i,j)^2 > \beta \mathbf{K}(i,i) \mathbf{K}(j,j)
\end{equation}
this choice of $\beta$ serves to disconnect nodes in the graph when the elements connecting them varies in stiffness by roughly 1 order of magnitude or more from any surrounding elements. Thus nodes in void regions are aggregated together, as are nodes in solid regions of the structure, but nodes are not connected across high gradients of element stiffness (i.e. along the edges of structural elements). We do not claim that this is necessarily an optimal choice of the parameter, merely that it is a very effective one in our experience.}
\delete{To construct the prolongation operator connecting the first two levels of multigrid, candidate vectors (the rigid body modes before boundary conditions are applied) are smoothed with 4 iterations of block Gauss-Seidel (1 block per node). No such candidate smoothing is applied to any of the subsequent prolongation operators connecting lower levels of the hierarchy. However,} Once assembled, each tentative prolongation operator is improved with a \delete{single pass of a} weighted Jacobi smoother. The system of linear equations for displacements is solved with \delete{preconditioned conjugate gradients}{\add GMRES \citep{Saad1986}} to a relative residual tolerance of 1e-7 in the cantilever problems and 1e-8 in the other two. The initial guess for the iterative solver is set to the displacements calculated in the last topology iteration.

For the optimization itself we use continuation on the penalty parameter, with the penalty initially set to 1 and gradually increased to 4. The number of iterations at each penalty value and the step size for the penalty values varies slightly between the different examples. Penalization is performed according to the modified SIMP method \citep{Sigmund1997DesignMethod}, with minimum stiffness set to 1e-10. The design variable updates are calculated using the Method of Moving Asymptotes (MMA) \citep{Svanberg1987}. In each case we use a density filter with radius equal to 1.5 times the element dimensions. \delete{The near constant coarse grid size and decreasing feature size will help to demonstrate how both multigrid approaches fare at representing finer features on coarse grids.}

\input{2DCantilever.tex}

\input{2DStability.tex}

\input{VaryingGrid.tex}

\input{3DCantilever.tex}

%% file: 2DCantilever.tex
\subsection{2D cantilever beam}
\label{section:2D_Cantilever}
\begin{figure*}
	\centering
	\hspace*{\fill}
	\begin{subfigure}[t]{0.45\textwidth}
    	\centering
        \includegraphics[width=0.922\linewidth]{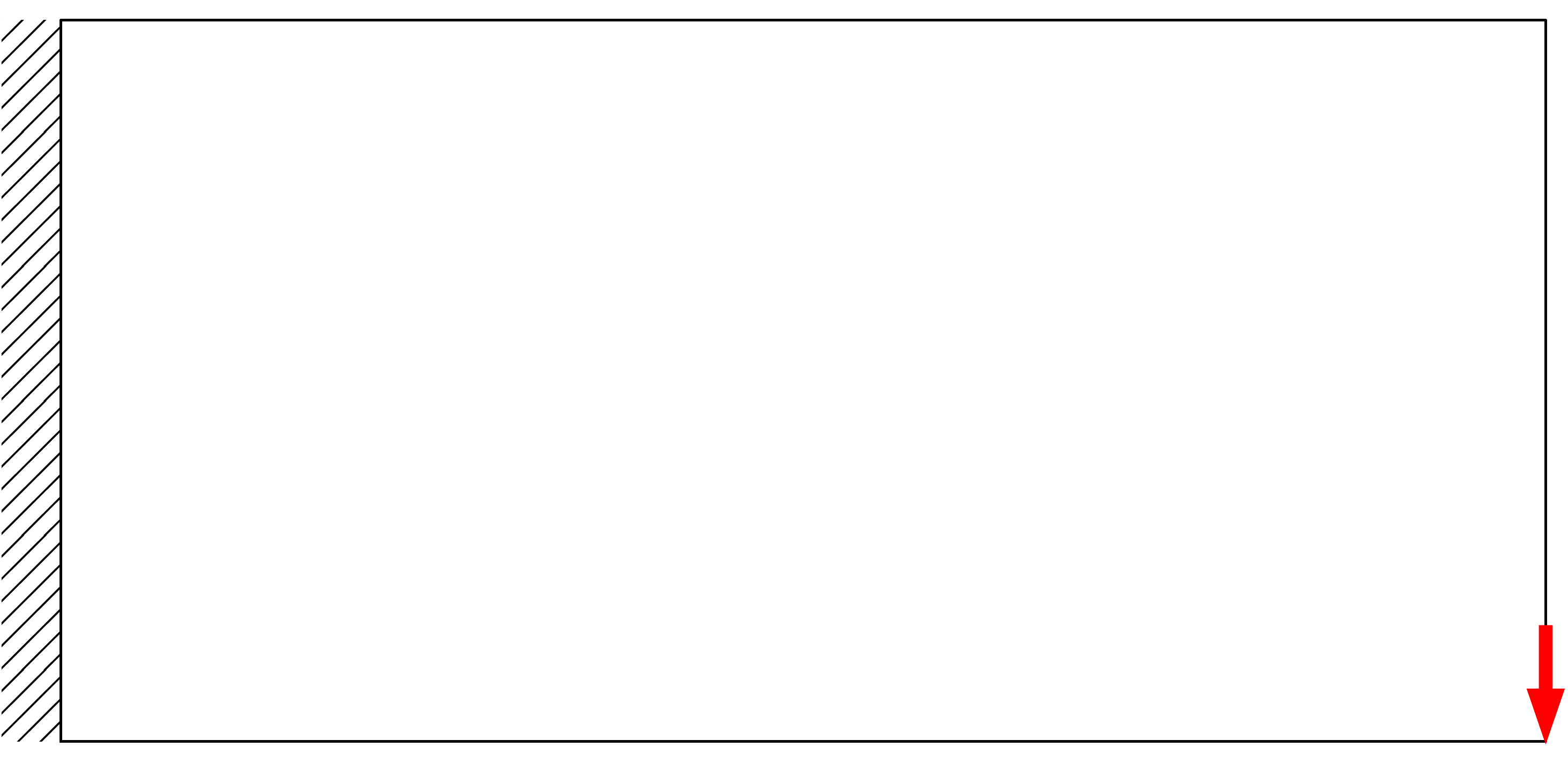}
	\end{subfigure}
    \hfill
	\begin{subfigure}[t]{0.45\textwidth}
        \centering
        \includegraphics[width=0.87\linewidth]{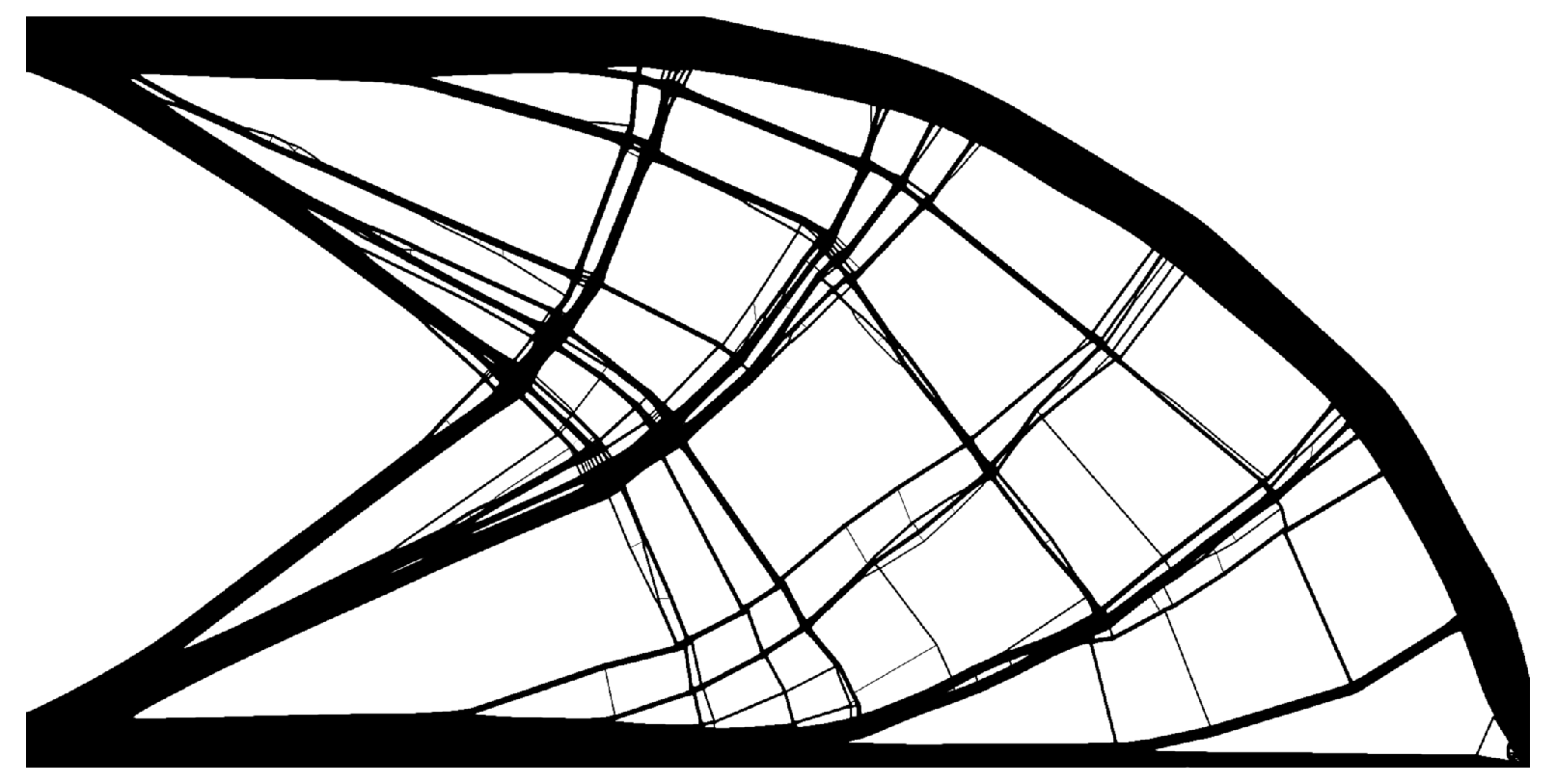}
	\end{subfigure}
	\hspace*{\fill}
    \caption{Design domain and optimized structure for 2D cantilever beam problem. The left edge of the domain is fixed and a load is applied at the middle of the right edge.}
    \label{fig:Cantilever}
\end{figure*}

The first example is compliance minimization for a cantilever domain in 2D with an aspect ratio of 2:1. The design domain and {\add resulting structure} \delete{result of the highest resolution optimization} are illustrated in Figure \ref{fig:Cantilever}. \delete{We run the problem at four different resolutions: 96x48 elements, 192x96 elements, 384x192 elements, and 768x384 elements.} {\add The domain is discretized at a resolution of 2048x1024 elements and the optimization is run in parallel across 128 processors.} We perform \delete{80} {\add 20} optimization iterations for each penalty value and the penalty is increased in increments of 0.25 {\add from 1 to 4. The maximum volume fraction is set to 0.4.}

\begin{figure}
	\centering
	\begin{subfigure}[t]{\columnwidth}
    	\centering
        \includegraphics[width=0.85\linewidth]{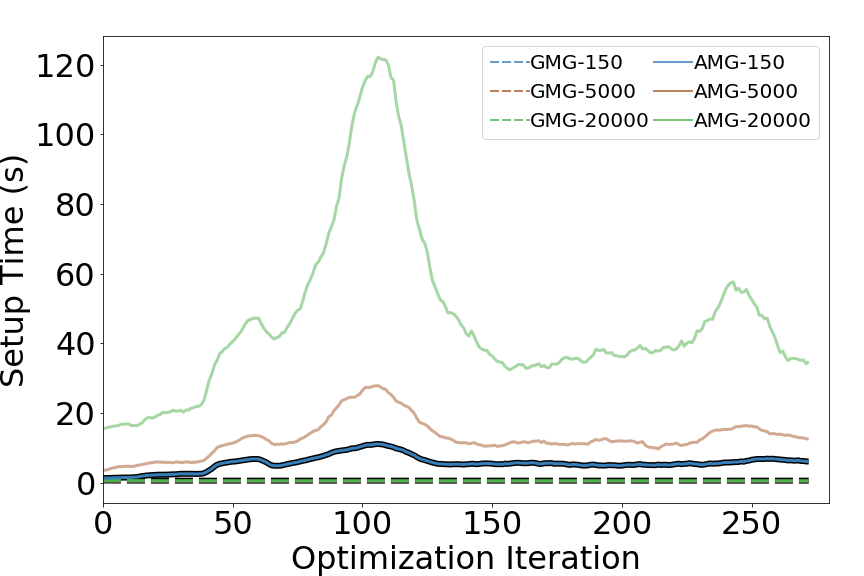}
        \caption{Time to set up the preconditioner.}
    \end{subfigure}

	\begin{subfigure}[t]{\columnwidth}
        \centering
        \includegraphics[width=0.85\linewidth]{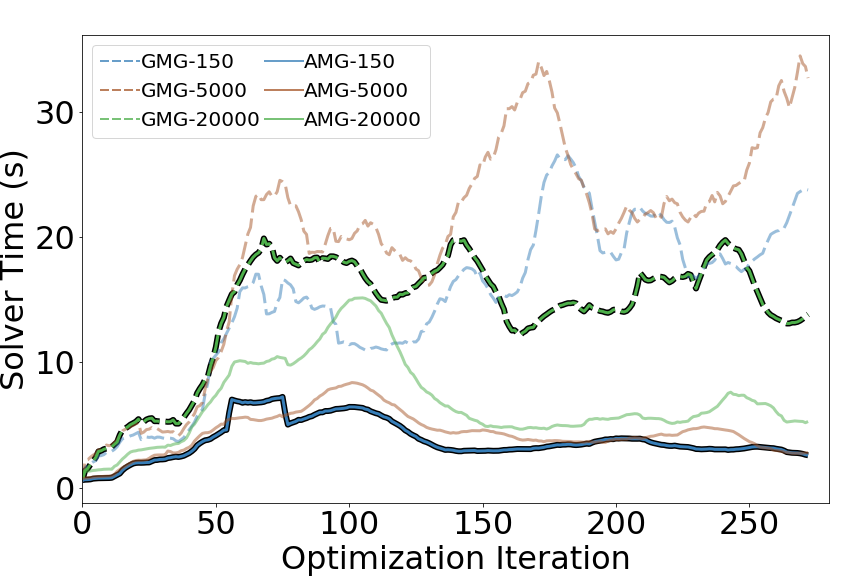}
        \caption{Time to solve the linear system for displacements.}
    \end{subfigure}
    
	\begin{subfigure}[t]{\columnwidth}
        \centering
        \includegraphics[width=0.85\linewidth]{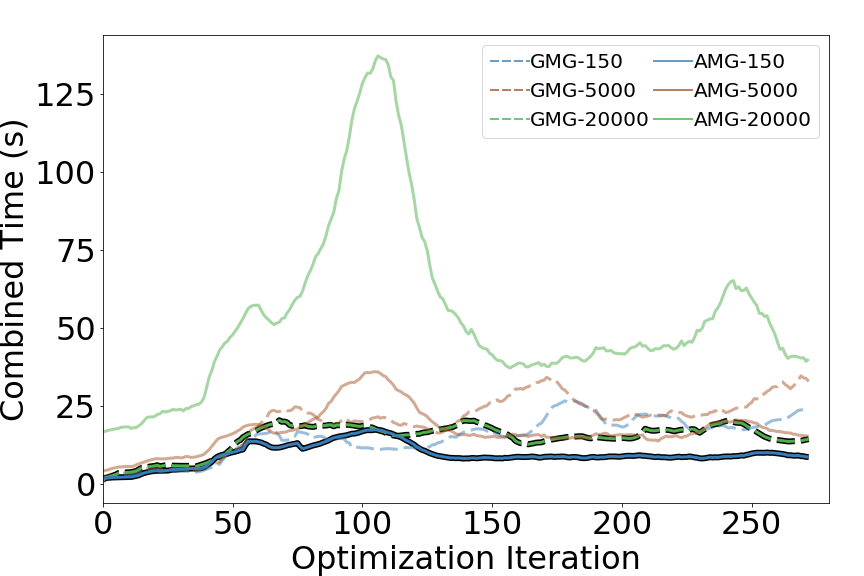}
        \caption{Combined time to set up and solve.}
    \end{subfigure}

    \begin{subfigure}[t]{\columnwidth}
        \centering
        \includegraphics[width=0.85\linewidth]{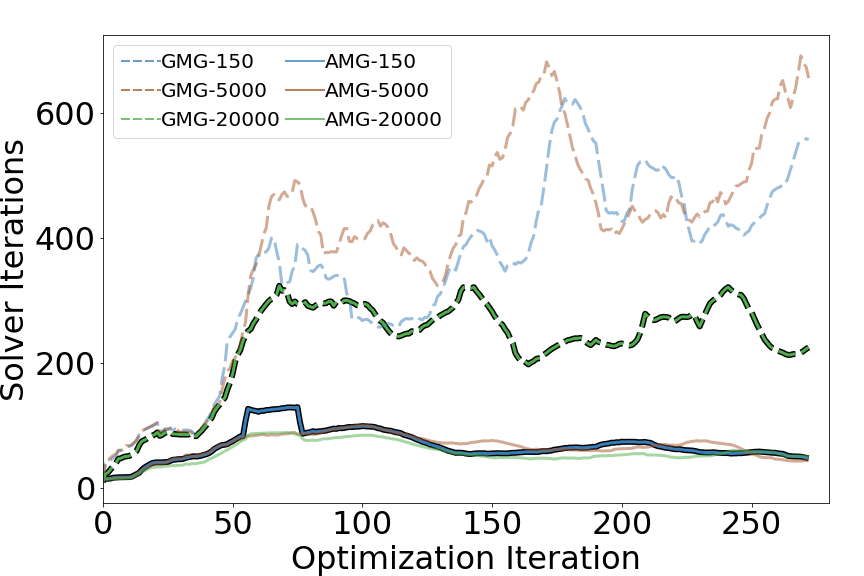}
        \caption{Iterations of the linear solver until convergence.}
    \end{subfigure}
    \caption{Performance of the preconditioners for the 2D cantilever problem with varying coarse grid size. AMG performance is displayed with solid lines and GMG with dashed lines. The best setup for AMG and GMG are emboldened for easier comparison.}
    \label{fig:2D_Cantilever_Coarse_Size}
\end{figure}

{\add We start by using weighted point-block jacobi (one block per node) with a weight of 0.5 as a smoother on all levels, and the coarse grid is solved with an LU decomposition on a single process. Figure \ref{fig:2D_Cantilever_Coarse_Size} demonstrates how the two multigrid versions perform as the coarse grid size is increased. The plots show a smoothed trendline rather than actual performance at each optimzation iteration for easier comparison between methods. We constrain each hierarchy to have coarse grid sizes less than 150, 5000, or 20000 dofs. This corresponds to 9, 6, or 5 levels in the GMG hierarchy and a similar (though varying) number for the AMG hierarchy.

As expected, for both methods the setup cost increases with increased size of the coarse grid, and the AMG solver is about 50 times more expensive to set up than GMG in each case. However, it is also important to note that as the coarse grid decreases in size, the GMG solver requires more iterations to solve. By the end of the optimization, the GMG preconditioner with the largest coarse grid requires several hundred iterations less than the ones with smaller coarse grids. In contrast, the AMG solver requires an approximately constant number of iterations, only varying by about 15 iterations from the coarsest to finest grid. For the GMG solver, the increased cost of additional iterations is enough to offset the decreased setup cost for smaller coarse grids. As a result, the AMG preconditioner is most efficient for the smallest coarse grid size, while the total cost of the GMG preconditioner is nearly independent of coarse grid size. These results suggest that the AMG preconditioner is most effective for a minimal coarse grid size, while the GMG preconditioner is most effective for a large one. Comparing the most effective choices for AMG vs. GMG (shown in bold), we see that the high setup cost of AMG is more than offset by the lower number of iterations to convergence. Overall, using the AMG preconditioner takes about 50\% less time to converge than using the GMG preconditioner.}

\begin{figure}
	\centering

	\begin{subfigure}[t]{\columnwidth}
        \centering
        \includegraphics[width=0.85\linewidth]{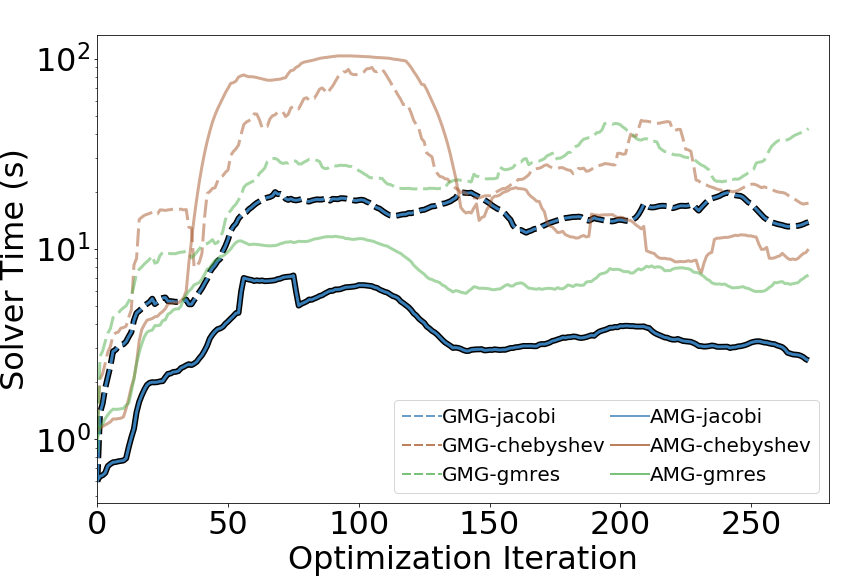}
    \end{subfigure}

    \begin{subfigure}[t]{\columnwidth}
        \centering
        \includegraphics[width=0.85\linewidth]{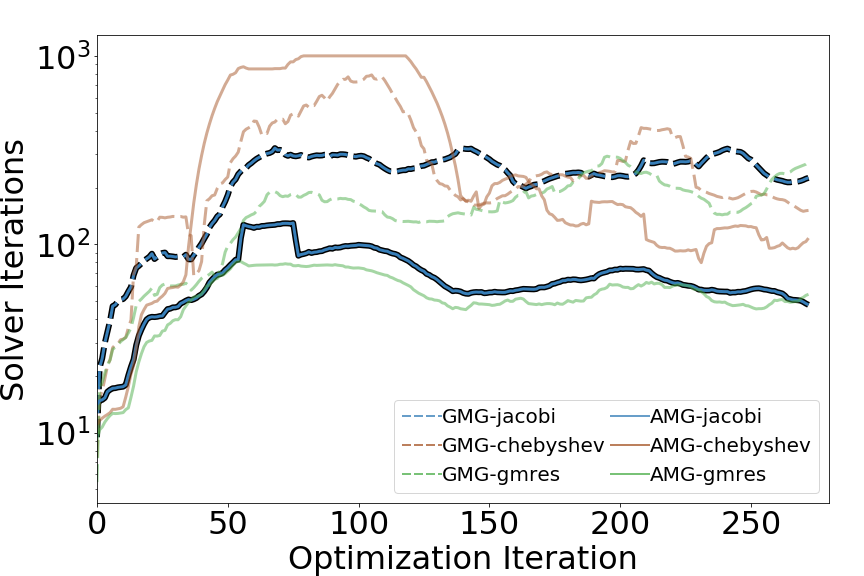}
    \end{subfigure}
    \caption{Performance of the preconditioners for the 2D cantilever problem with different smoothers. AMG performance is displayed with solid lines and GMG with dashed lines. The best setup for AMG and GMG are emboldened for easier comparison.}
    \label{fig:2D_Cantilever_Smoothers}
\end{figure}

{\add We also use this opportunity to compare the performance of three different smoothing techniques, SOR-preconditioned chebyshev iterations (the default for PETSC MG preconditioners, used in \citep{Aage2017Giga-voxelDesign}), SOR-preconditioned GMRES iterations (used in \citep{Aage2015TopologyFramework}), and weighted jacobi (used in \citep{Amir2014OnOptimization}). To accommodate the GMRES smoother we use the flexible GMRES method \citep{Saad1993AAlgorithm} as a linear solver. Given that difference in setup cost for these smoothers is negligible, we present the change in solver time and number of iterations for the most effective AMG and GMG setup for this problem in Figure \ref{fig:2D_Cantilever_Smoothers}. For both the GMG and AMG hierarchy, the GMRES smoother requires the least number of iterations while the jacobi smoother results in the least time to solve. 

Interestingly, the chebyshev smoother sees a huge spike in the number of required iterations from roughly steps 40 to 130 of the optimization, and for the AMG hierarchy often fails to converge within 1000 iterations. Though not explored here, a more robust convergence criteria, such as the one used in \citep{Amir2014OnOptimization}, or better tuning of the smoother parameters may reduce the number of iterations in these cases. Performance could also be improved with better estimates of the operator eigenvalues (which are calculated internally by the PETSC framework); however techniques to improve the estimate are beyond the scope of this paper.}

%% file: 2DStability.tex
\subsection{{\add 2D} Column Stability}
\label{section:2D_Column}
\begin{figure}
	\centering
	\hspace*{\fill}
	\begin{subfigure}[c]{0.45\columnwidth}
		\centering
		\includegraphics[width=0.6\linewidth]{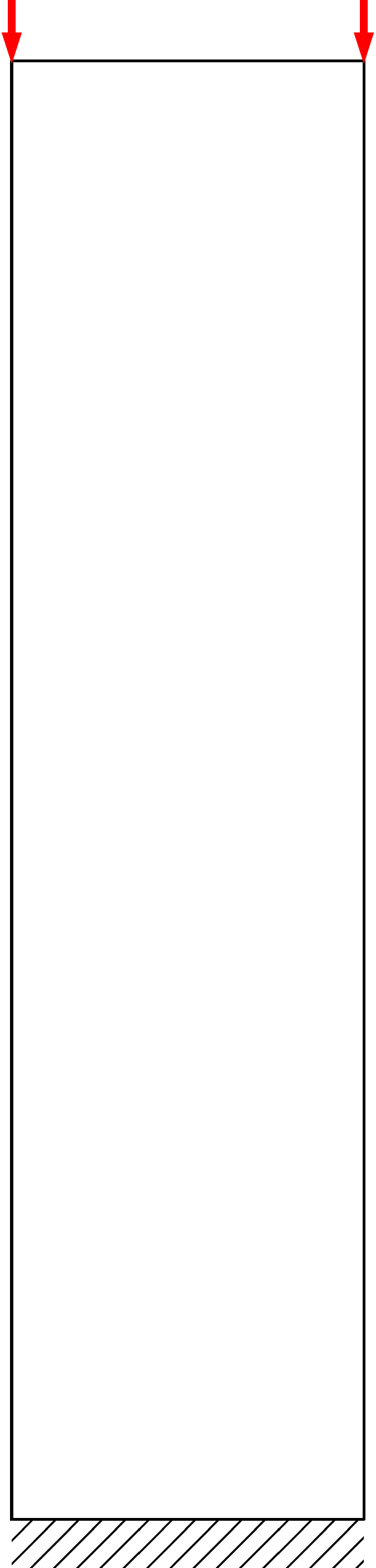}
	\end{subfigure}
    \hfill
	\begin{subfigure}[c]{0.45\columnwidth}
		\centering
		\includegraphics[width=0.62\linewidth]{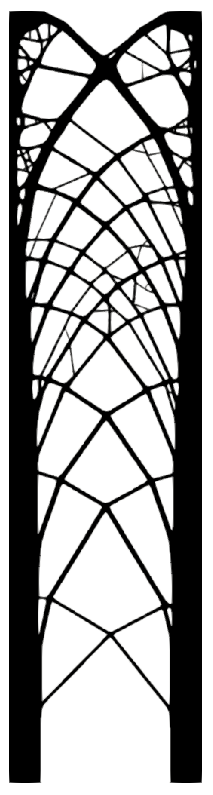}
	\end{subfigure}
	\hspace*{\fill}
    \caption{Column stability problem domain and result.}
    \label{fig:Stability}
\end{figure}

The next problem is less common in the literature, but will help us to further differentiate the performance of the two multigrid approaches. Here we perform stability optimization of a column as described in \citep{Bendse2003}. The design domain has aspect ratio 4:1 and we run the problem at \delete{three different resolutions: 32x128 elements, 64x256 elements, and 96x392 elements} {\add a resolution of 256x1024 elements across 64 processors}. We perform \delete{70} {\add 30} optimization iterations for each penalty value and the penalty is increased in increments of \delete{0.25}{\add 0.125 to a value of 4. As this is not enough for the structure to fully converge, we continue to increase the penalty by a value of 0.25 every 40 iterations up to a final penalty of 12.}. {\add The volume fraction for this example is again set to 0.4.} The design domain and boundary conditions, as well as a sample optimized shape\delete{ from the largest problem}, are shown in Figure \ref{fig:Stability}.

When optimizing for stability we have to solve both an eigenvalue problem for the structural performance and an adjoint problem for the sensitivities. {\add Other works have successfully made use of GMG preconditioning for adjoint problems \citep{Alexandersen2016LargeConvection}, though they do not compare the performance to that of an AMG preconditioner as we will here.} Whereas compliance optimization requires only a single solution to a linear system each time the preconditioner is constructed, stability optimization requires multiple solutions to the same linear system with different right-hand-sides. Thus, we use this example to demonstrate how the preconditioners perform when the improvement to the conditioning of the system is more important relative to the setup cost. {\add Some approaches, such as Krylov subspace recycling \citep{Parks2006RecyclingSystems,Wang2007Large-scaleRecycling}, have been established to reduce the cost of solving the same linear system with multiple right-hand-sides, though they aren't applicable when solving for eigenvalues (easily the most expensive step) and the potential benefit for the adjoint solutions is not considered in this paper.}

When solving for displacements we use the same procedure as before, however the procedure to solve the adjoint problem is slightly different. We retain the same iterative solver and preconditioners, but now we use an all-zero initial condition \delete{and the residual tolerance is also relaxed to 1e-5.}{\add as it is difficult to consistently provide a good initial guess due to the potential for modes switching. This also gives us the opportunity to compare how the preconditioners perform when the initial guess lacks intuition from previous results.} To solve the eigenvalue problem we use the generalized Davidson method with the same preconditioners as before. At every optimization step we calculate six eigenvalues \delete{to a relative residual tolerance of 1e-6 for a maximum of 1e3 eigenvalue iterations} {\add and consider them converged when the relative residual falls below 1e-6 or the relative change in the approximated eigenvalue between consecutive iterations is less than 1e-13.}

\delete{The performance when solving for displacements is shown in Figure \ref{fig:Stability_Performance}. The trends are again consistent with what we identified in the 2D cantilever example, namely that the AMG solvers have a higher setup cost, but overall the time to solve for displacements is lower thanks to a significant reduction in the number of CG iterations. The main difference for the stability optimization problem is that the AMG preconditioner performs noticeably worse than the block AMG preconditioner. This is due to the fact that the block smoothing allows the preconditioner to more accurately solve for the displacements of supernodes on the coarse grids. Quantitatively, the time to setup the preconditioner and solve for displacements is consistently twice as much for GMG than block AMG and 50\% higher for AMG than block AMG.}

\begin{figure}
	\centering
	\begin{subfigure}[t]{\columnwidth}
    	\centering
        \includegraphics[width=0.85\linewidth]{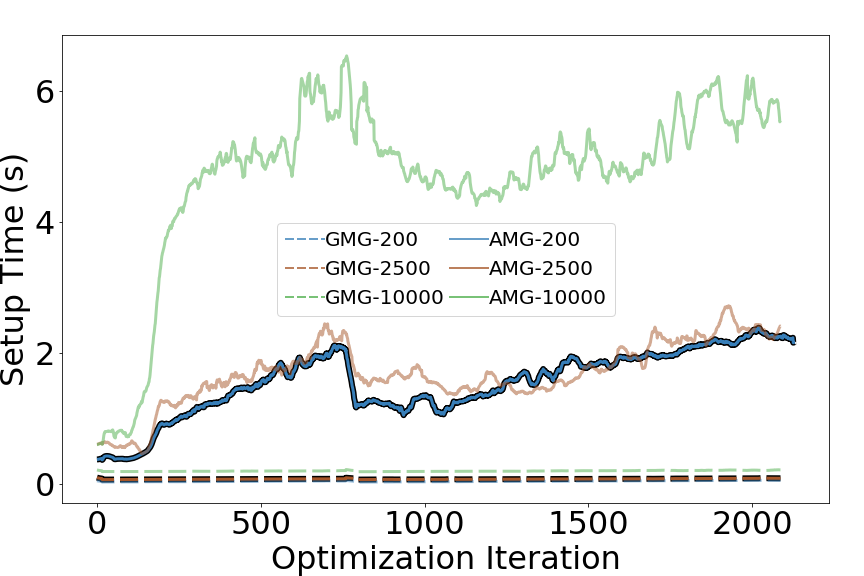}
        \caption{Time to set up the preconditioner.}
    \end{subfigure}

	\begin{subfigure}[t]{\columnwidth}
        \centering
        \includegraphics[width=0.85\linewidth]{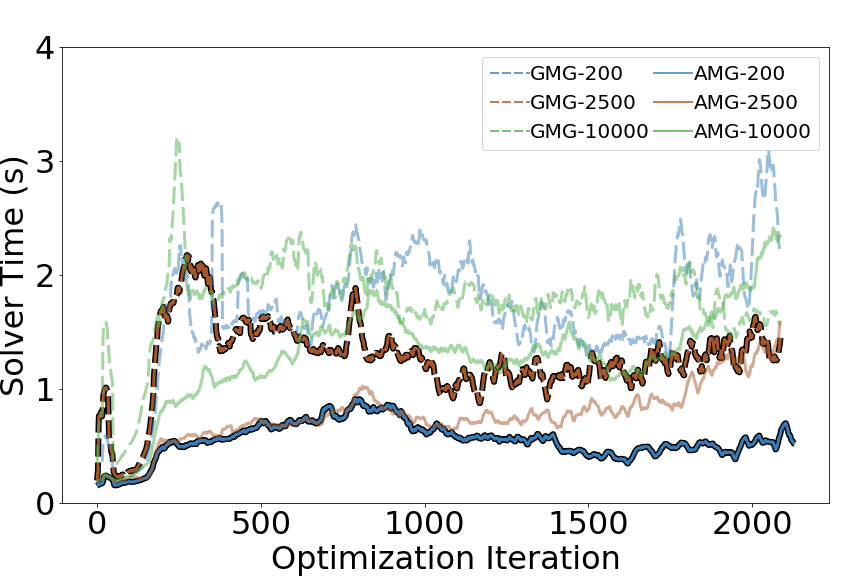}
        \caption{Time to solve the linear system for displacements.}
    \end{subfigure}
    
	\begin{subfigure}[t]{\columnwidth}
        \centering
        \includegraphics[width=0.85\linewidth]{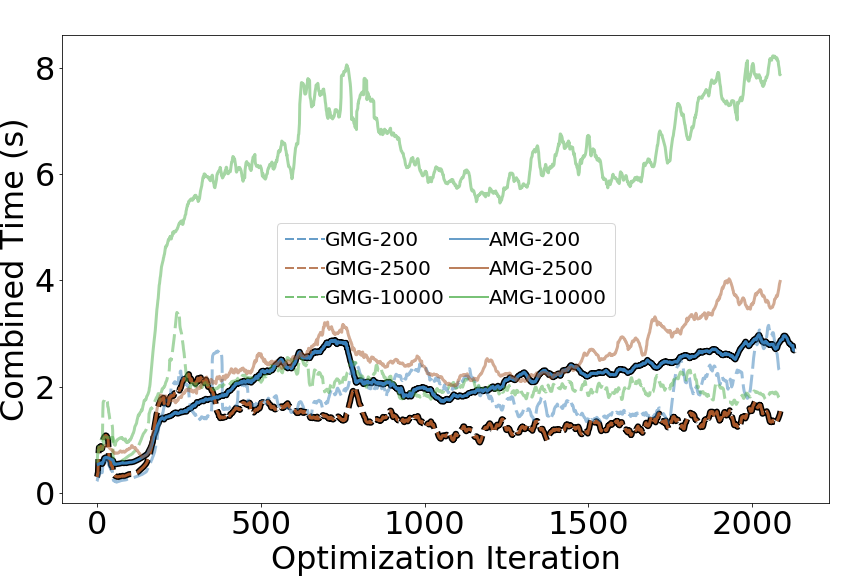}
        \caption{Combined time to set up and solve.}
    \end{subfigure}

    \begin{subfigure}[t]{\columnwidth}
        \centering
        \includegraphics[width=0.85\linewidth]{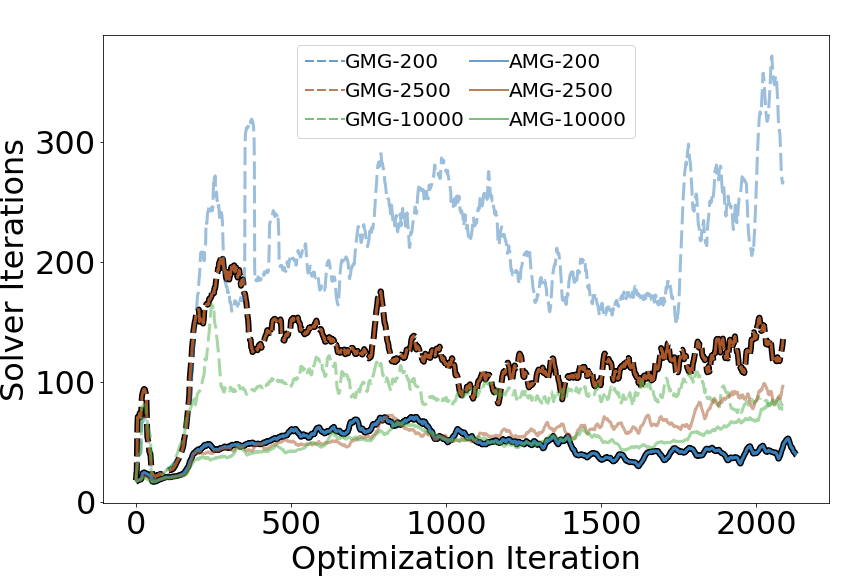}
        \caption{Iterations of the linear solver until convergence.}
    \end{subfigure}
    \caption{Performance of the preconditioners {\addd when solving for displacements in}\deleted{for} the 2D stability problem with varying coarse grid size. AMG performance is displayed with solid lines and GMG with dashed lines. The best setup for AMG and GMG are emboldened for easier comparison.}
    \label{fig:2D_Stability_Coarse_Size}
\end{figure}

{\add We compare the performance of each method when solving for displacements using the weighted jacobi smoother with varying coarse grid size in Figure \ref{fig:2D_Stability_Coarse_Size}. As the fine scale problem is slightly smaller than before, we compare GMG preconditioners with 4, 5, or 7 levels (8514, 2210, and 170 dofs, respectively) to AMG preconditoners restricted to 10000, 2500, and 200 dofs on the coarse grids.

The overall trends for the AMG smoother are the same as for the cantilever problem, iteration counts are very consistent across grid sizes, but the overall time is smallest for the smallest coarse grid size. For the GMG solver the trend changes slightly. While iteration counts still decrease with increasing coarse grid size, the number of iterations for the smallest grid sizes is only about half as much as in the cantilever problem. As a result, there is less room for improvement from increasing the coarse grid size and the tradeoff is not enough to make the 4-level operator the most effective. Instead the 5-level operator is the most efficient overall, though the 7-level solver provides very similar compute time because the increase in iterations from additional refinement is modest as well.}

\begin{figure}
	\centering
	\begin{subfigure}[t]{\columnwidth}
    	\centering
        \includegraphics[width=0.85\linewidth]{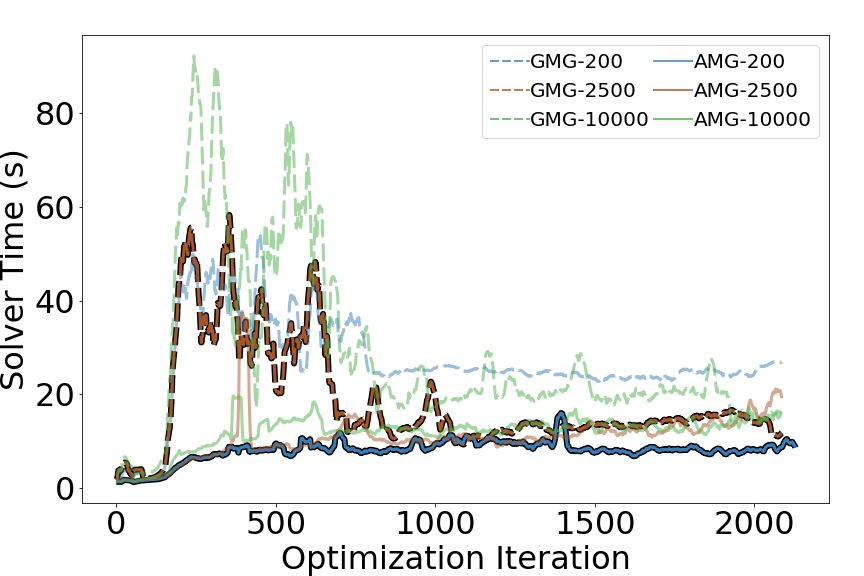}
        \caption{Time to solve for eigenvalues.}
    \end{subfigure}

	\begin{subfigure}[t]{\columnwidth}
        \centering
        \includegraphics[width=0.85\linewidth]{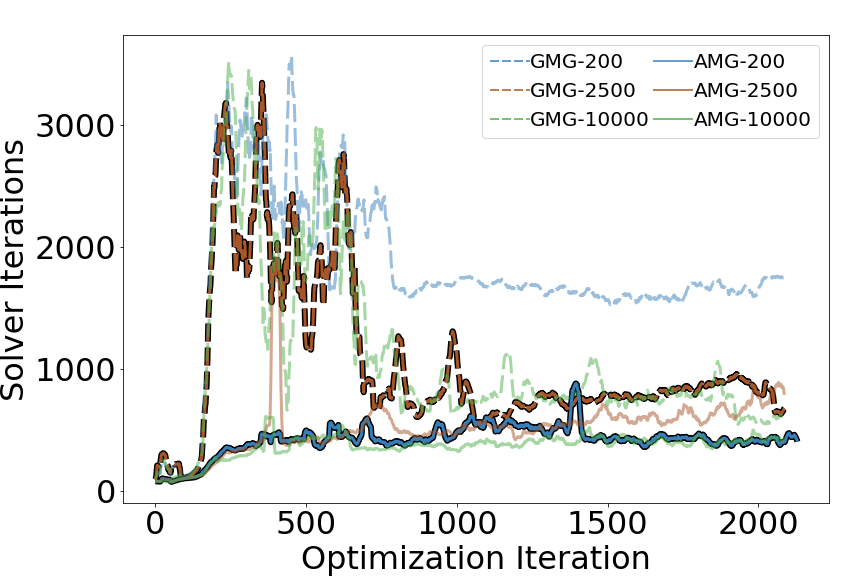}
        \caption{Number of iterations to solve for eigenvalues.}
    \end{subfigure}
    \caption{Performance of the preconditioners for the 2D stability problem with varying coarse grid size when solving for eigenvalues. AMG performance is displayed with solid lines and GMG with dashed lines. The best setup for AMG and GMG are emboldened for easier comparison.}
    \label{fig:2D_Eigensolver_Coarse_Size}
\end{figure}

\begin{figure}
	\centering
	\begin{subfigure}[t]{\columnwidth}
    	\centering
        \includegraphics[width=0.85\linewidth]{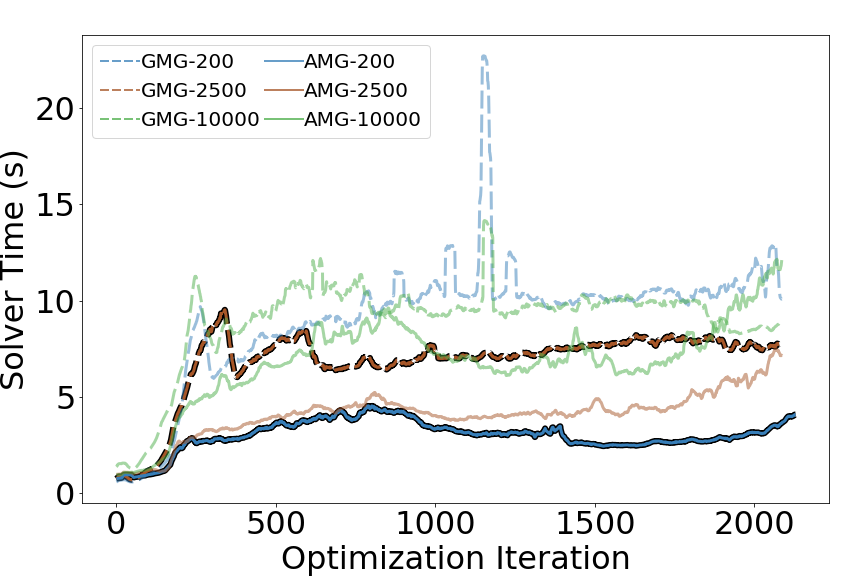}
        \caption{Time to solve the adjoint problems.}
    \end{subfigure}

	\begin{subfigure}[t]{\columnwidth}
        \centering
        \includegraphics[width=0.85\linewidth]{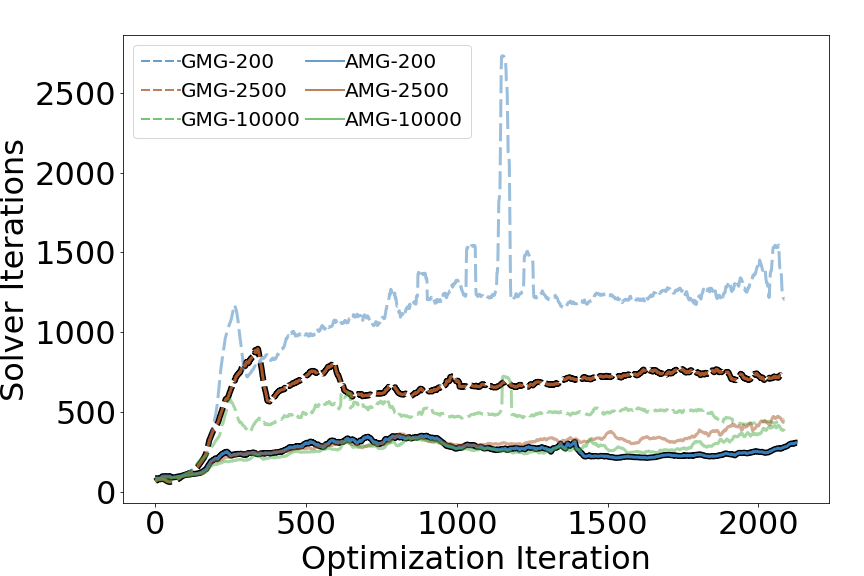}
        \caption{Number of iterations to solve the adjoint problems.}
    \end{subfigure}
    \caption{Performance of the preconditioners for the 2D stability problem with varying coarse grid size when solving for the adjoint problems. AMG performance is displayed with solid lines and GMG with dashed lines. The best setup for AMG and GMG are emboldened for easier comparison.}
    \label{fig:2D_Adjoint_Coarse_Size}
\end{figure}

{\add We similarly compare the performance of the preconditioners with varying coarse grid size when solving for eigenvalues in Figure \ref{fig:2D_Eigensolver_Coarse_Size} and when solving the adjoint problems in Figure \ref{fig:2D_Adjoint_Coarse_Size}. The relative trends are the same as when solving for displacements and the same setups that were most effective for displacements are again most effective here, but we note intermediate spikes when solving for eigenvalues using the AMG preconditioner. These spikes occur for the same reasons as discussed for the Chebyshev smoother previously (though we are using weighted Jacobi in this case), stagnation of the eigensolver for a particularly challenging intermediate topology. However, even with the spike in number of iterations, it still arrives at a solution in fewer iterations than any of the GMG preconditioners at that stage.}

\begin{figure}
	\centering
    \includegraphics[width=0.85\linewidth]{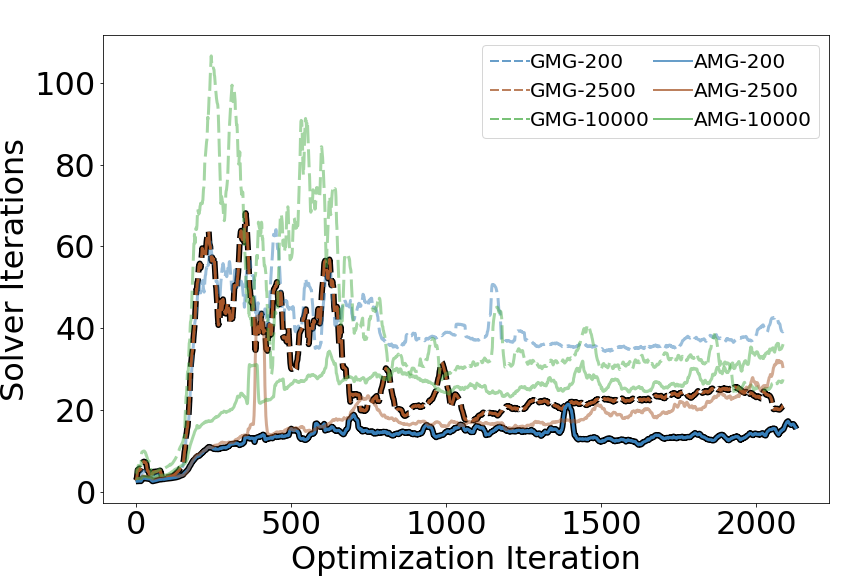}
    \caption{Total time to set up the preconditioners and solve for displacements, eigenmodes, and adjoint solutions in the 2D stability problem with varying coarse grid size. AMG performance is displayed with solid lines and GMG with dashed lines. The best setup for AMG and GMG are emboldened for easier comparison.}
    \label{fig:2D_Overall_Coarse_Size}
\end{figure}

{\add Comparing the most effective version of the AMG and GMG preconditioners, The high setup cost for AMG is again offset by the increased iterations for the GMG preconditioner. For this problem though, the number of iterations to solve for displacements remains low enough that GMG is most effective to setup and solve a single linear system. However, the much more expensive steps to solve for eigenvalues and the associated adjoint problems make the AMG preconditioner more effective overall by drastically reducing the cost of these steps (Figure \ref{fig:2D_Overall_Coarse_Size}). While the AMG preconditioner takes almost twice as long to setup and solve a single linear system by the end, the difference is only about 1.5 seconds. In contrast, the difference in times to calculate the eigenvalues alone can exceed 40 seconds.

Looking at the performance when solving the adjoint problems, we also note that regardless of the preconditioner used, the time and number of iterations to solve all 6 adjoint problems is only about 5 times more than solving the single system for displacements. This indicates that our choice to use an all-zero initial guess for these problems does not substantially harm performance. Though not shown here, it is the author's experience that using the previous eigenvectors as an initial searchspace for the eigensolver still leads to a major reduction in the number of iterations and should not be neglected.}

%% file: VaryingGrid.tex
\subsection{Varying Spaced Grid}
\label{section:Grid}
\begin{figure*}
    \centering
    
	\begin{subfigure}[t]{0.48\linewidth}
    	\hspace*{\fill}
    	\begin{subfigure}[t]{0.32\linewidth}
        	\centering
            \includegraphics[width=\linewidth]{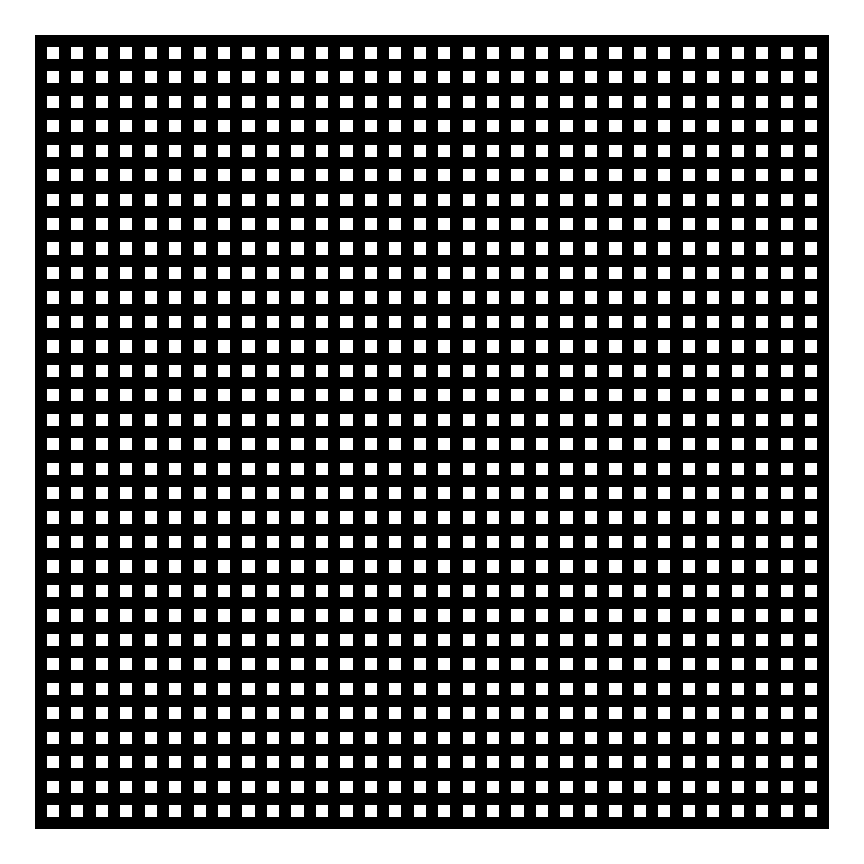}
        \end{subfigure}
        \hfill
    	\begin{subfigure}[t]{0.32\linewidth}
        	\centering
            \includegraphics[width=\linewidth]{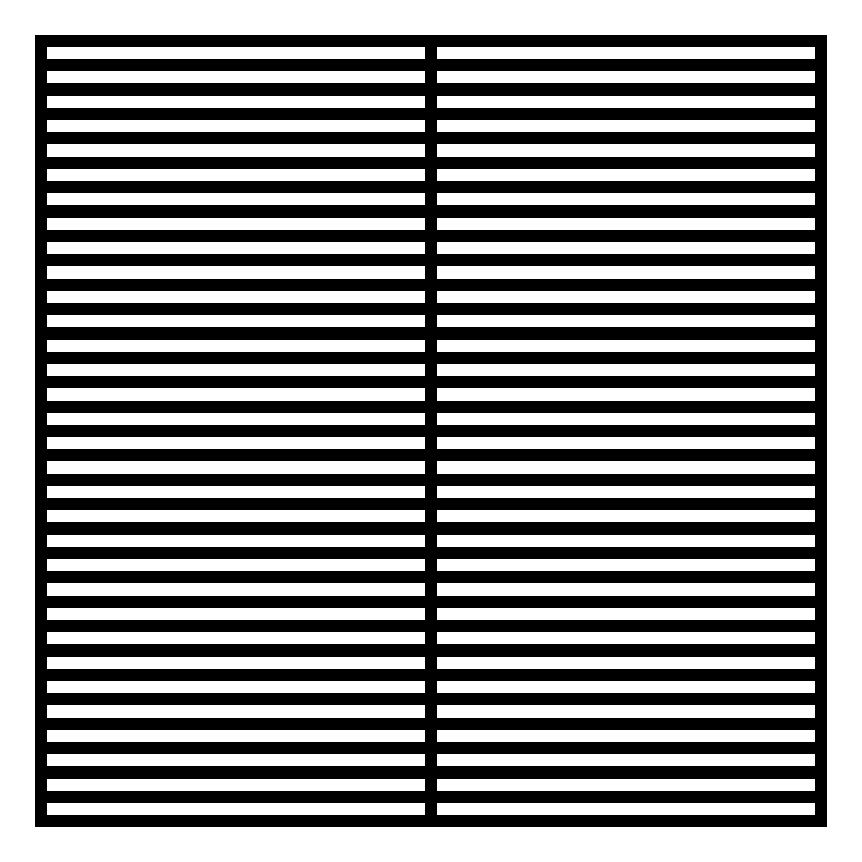}
        \end{subfigure}
        \hfill
    	\begin{subfigure}[t]{0.32\linewidth}
        	\centering
            \includegraphics[width=\linewidth]{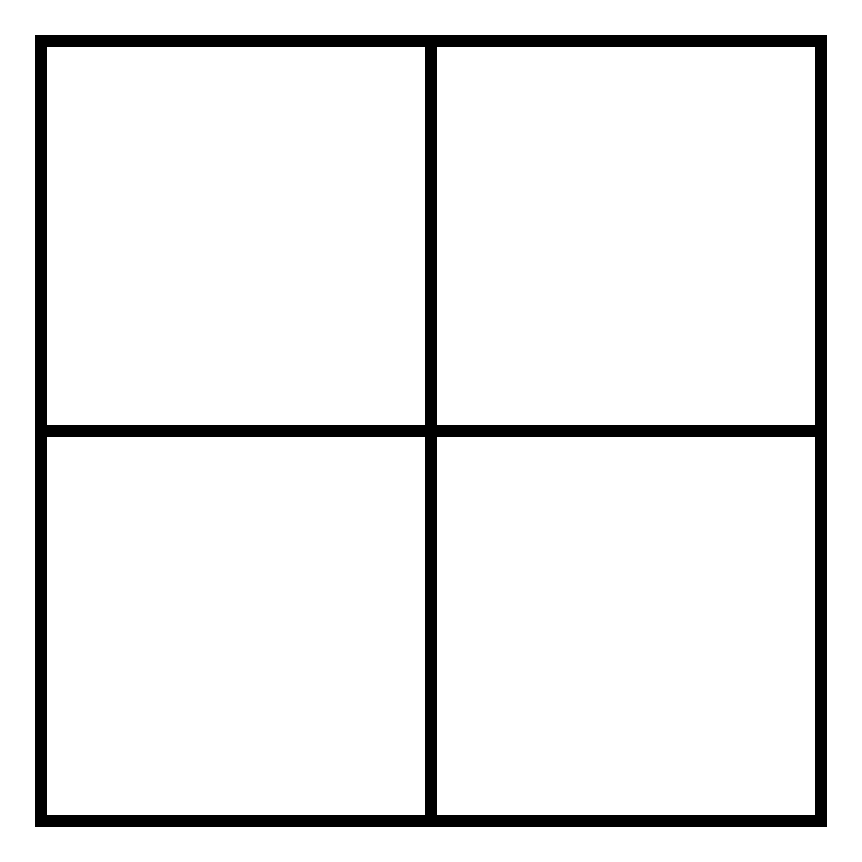}
        \end{subfigure}
    	\hspace*{\fill}
    	
        \caption{Sample grid structures. The first shows minimum spacing of columns and beams, and the last shows the maximum spacing. The middle figure shows maximum column spacing and minimum beam spacing.}
        \label{subfig:Grid}
    \end{subfigure}
    \hfill
	\begin{subfigure}[t]{0.48\linewidth}
        \centering
    	\hspace*{\fill}
    	\begin{subfigure}[t]{0.32\linewidth}
        	\centering
            \includegraphics[width=\linewidth]{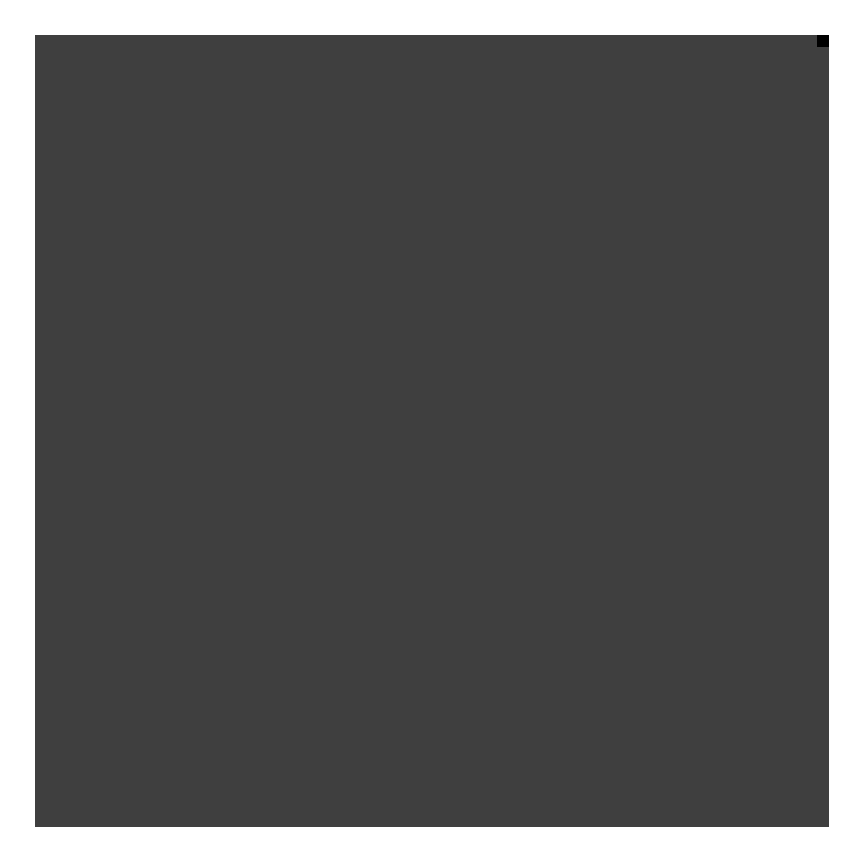}
        \end{subfigure}
        \hfill
    	\begin{subfigure}[t]{0.32\linewidth}
        	\centering
            \includegraphics[width=\linewidth]{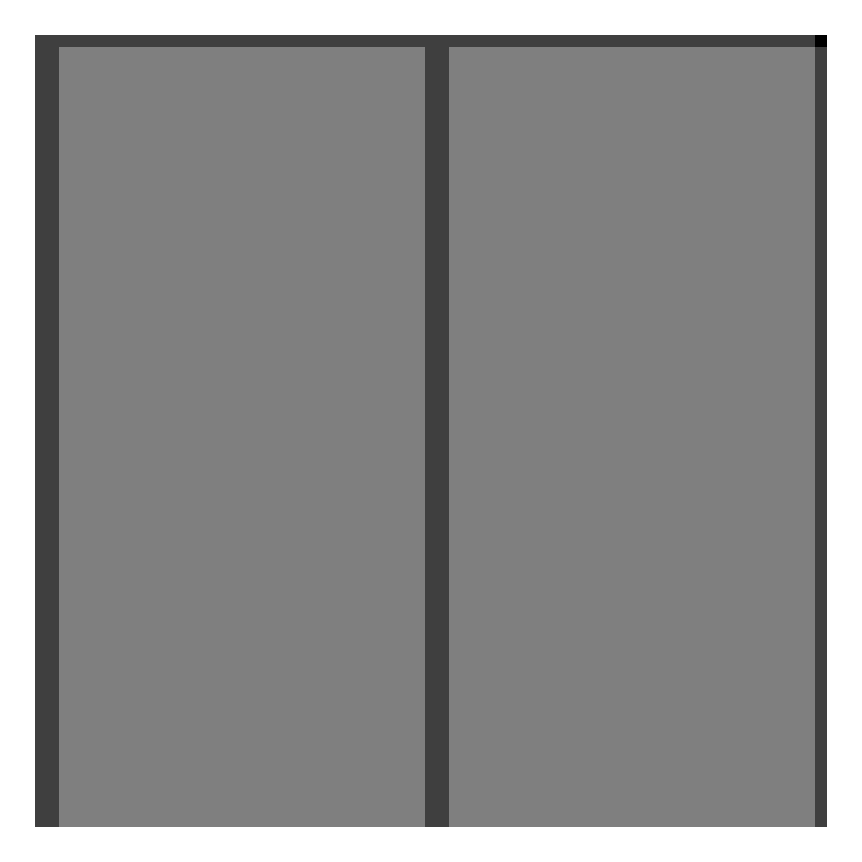}
        \end{subfigure}
        \hfill
    	\begin{subfigure}[t]{0.32\linewidth}
        	\centering
            \includegraphics[width=\linewidth]{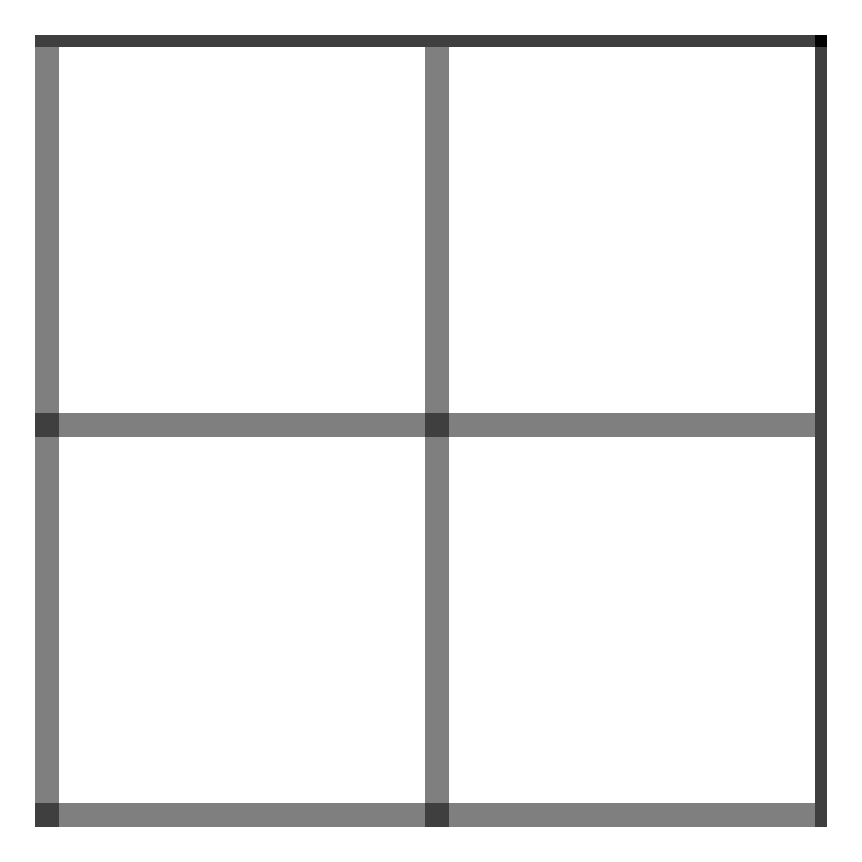}
        \end{subfigure}
    	\hspace*{\fill}
    	
        \caption{Sample grid structures as approximated on the coarse GMG grids.}
        \label{subfig:Coarsened_Grid}
    \end{subfigure}
    
	\begin{subfigure}[t]{\linewidth}
    	\begin{subfigure}[t]{0.48\linewidth}
            \centering
            \includegraphics[width=0.85\linewidth]{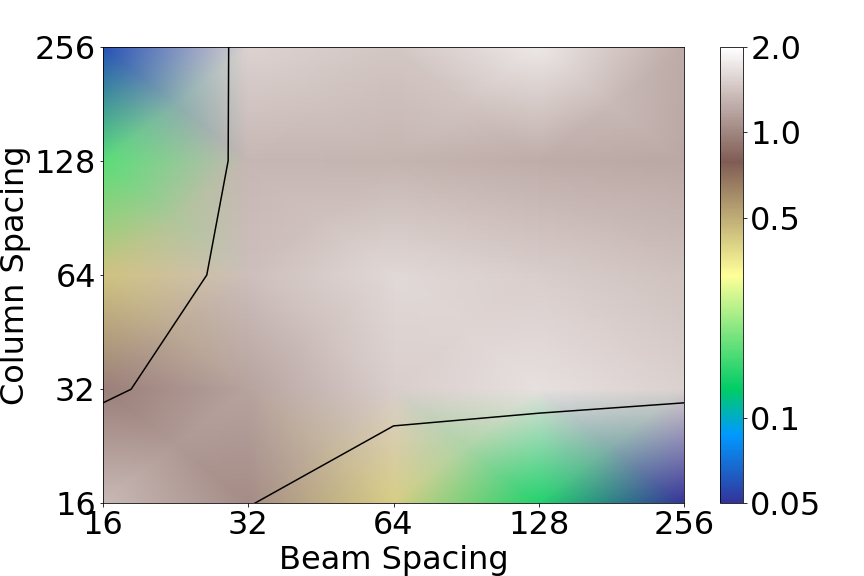}
        \end{subfigure}
        \hfill
        \begin{subfigure}[t]{0.48\linewidth}
            \centering
            \includegraphics[width=0.85\linewidth]{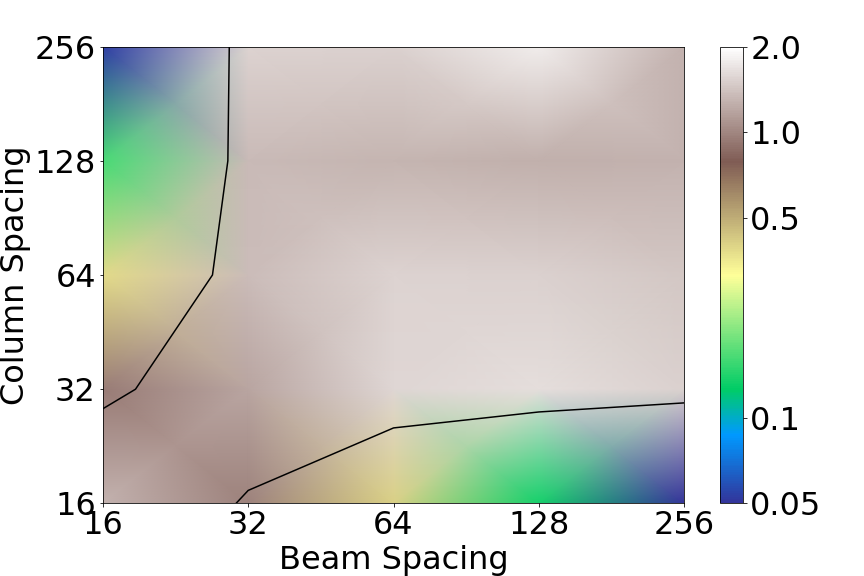}
        \end{subfigure}
        \caption{Ratios of cost of AMG compared to GMG for the 2D Grid. High values indicate GMG is more efficient, low values indicate AMG is more efficient.}
        \label{subfig:2D_Grid_BA-BG}
    \end{subfigure}

	\begin{subfigure}[t]{\linewidth}
        \begin{subfigure}[t]{0.48\linewidth}
            \centering
            \includegraphics[width=0.85\linewidth]{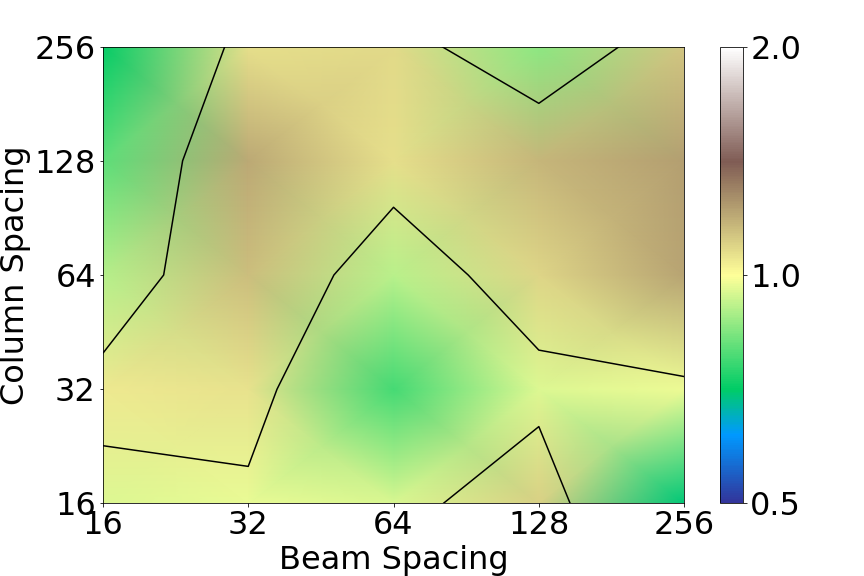}
        \end{subfigure}
        \hfill
        \begin{subfigure}[t]{0.48\linewidth}
            \centering
            \includegraphics[width=0.85\linewidth]{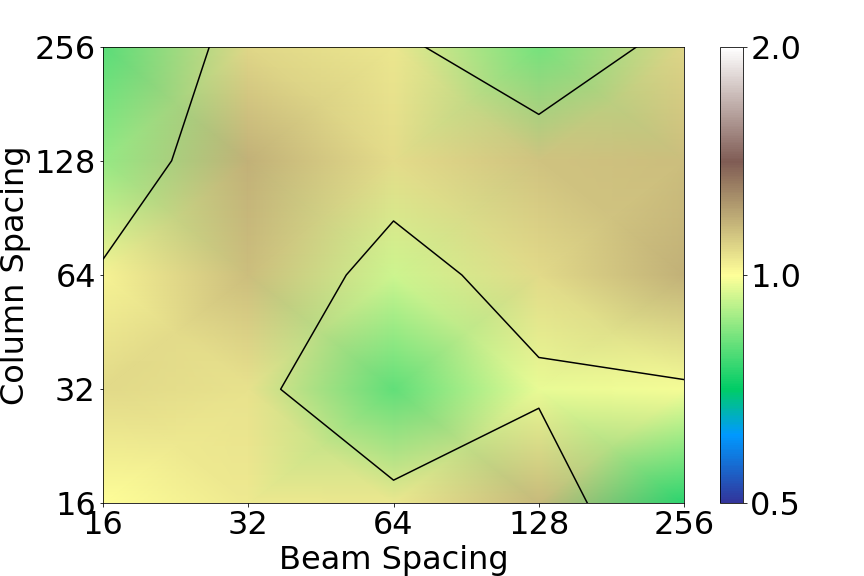}
        \end{subfigure}
        \caption{Ratios of cost of hybrid AMG compared to AMG for the 2D Grid. High values indicate AMG is more efficient, low values indicate the hybrid method is more efficient.}
        \label{subfig:2D_Grid_H-BA}
    \end{subfigure}
    \caption{Varying spaced grid example. Figures \ref{subfig:Grid} and \ref{subfig:Coarsened_Grid} show how the structures are represented on the finest and coarsest grids, respectively. Figures \ref{subfig:2D_Grid_BA-BG} and \ref{subfig:2D_Grid_H-BA} compare the performance of AMG, GMG a hybrid GMG-AMG method. The ratio of time spent in the linear solver is on the left and the ratio of solver iterations is shown on the right. The black line separates regions where each method outperforms the other.}
\end{figure*}

The final {\add 2D} problem is somewhat contrived, but serves to illustrate the exact situations where AMG outperforms GMG for topology optimization. Here, we use {\add a square domain of unit dimensions, and a fixed mesh resolution of 520x520 elements run on a single processor. We apply a uniform distributed load to the top of the domain, and fix all the degrees of freedom along the bottom. Rather than optimize the structure, we simply place 8-element-wide beams and columns at varying spacings of 16, 32, 64, 128, or 256 elements. A subset of the structures are shown in Figure \ref{subfig:Grid}. We construct an AMG preconditioner restricted to less than 1500 dofs on the coarse grid and a geometric hierarchy with 5 levels (2312 dofs on the coarse grid), and solve for displacements once. The corase grid is factorized with an LU decomposition and Weighted Jacobi smoothing is used on the remaining levels. Using 5 levels in the GMG preconditioner means that "elements" on the coarse grid have dimensions 16 times greater than the finest level. At this point the coarse representation is no longer a simple rediscretization of the fine grid, but instead performs some homogenization of material stiffness (Figure \ref{subfig:Coarsened_Grid}). Note that for the structures with the smallest feature spacing no void regions appear in the coarse representation, instead the maximum and minimum stiffness in the structure are only separated by a factor of 2 (as opposed to the void stiffness being 10 orders of magnitude smaller than solid material stiffness in the original representation).}
\delete{the same domain from the previous 2D cantilever example with aspect ratio of 2:1. However, now we fix the top 5\% of the domain to be a solid structure and apply a uniform distributed load to the top edge. Fixed supports of width 0.08 are applied with various uniform spacings to the bottom of the domain. In the case of just 2 supports at either end, we approximate the bridge problem described in \citep{Zegard2016BridgingManufacturing}, albeit only in 2D. As the number of supports increases, we see longer and thinner columns spanning from the top of the domain all the way to the supports at the bottom (Figure \ref{fig:Table_Results}). These columns and their associated low-energy bending modes prove troublesome for the GMG preconditioner, while the AMG preconditioner handles them without trouble. We run the optimization with 2, 4, 8, 16, and 25 supports. The case of 25 supports is equivalent to fixing every node on the bottom of the domain. The optimized structures from each case are shown in Figure \ref{fig:Table_Results}.}

{\add We compare the cost of the AMG and GMG preconditioners in Figure \ref{subfig:2D_Grid_BA-BG}. We only display the number of iterations and total time to set up the precondtioner and solve the system as the setup cost is nearly independent of the structure in this case and only depends on the method used (about 2.5 seconds for AMG and 0.5 seconds for GMG). As the figures show, the GMG performance is generally slightly better than the AMG performance in terms of both compute time and number of iterations, except when either the columns or beams are minimally spaced. When one set of features is minimally spaced and the other spacing increases, GMG quickly becomes very inefficient compared to AMG. As one spacing increases relative to the other, the features become thin and elongated, developing very smooth deformation modes. These modes are difficult to remove with smoothers, and the features must be coarsened repeatedly before they can be removed from the approximate solution.

However, when one spacing is set to the minimum, all the features are effectively merged at the 5th level of geometric refinement (Figure \ref{subfig:Coarsened_Grid}), and the coarse grid representation is nearly the same regardless of the original structure. Whereas sitffnesses vary by several orders of magnitude on the finer levels, they only vary by a factor of 2 on the coarse grid. If the spacing between features is even thinner, merging of features can happen with even fewer levels of refinement. Similarly, adding more levels to the GMG preconditioner will degrade performance for the cases where features are spaced even farther apart. In these scenarios the AMG preconditioner is much more effective as it naturally avoids merging features when coarsening. This concept also explains why the GMG preconditioner performed relatively better in the column example than it did in the cantilever, as optimizing for stability inhibits the development of long, thin features. Similarly, using a large GMG coarse grid size, which is less prone to prematurely merging features, was more beneficial in the cantilever problem than the stability one.

The phenomenon described here also implies the merit of a hybrid AMG approach similar to the one used in \citep{Aage2017Giga-voxelDesign}. For our grid example the structural features are all 8 elements wide and the minimum dimension of voids is at least that, meaning that 3 levels of geometric refinement can be performed without any risk of merging features. Given that the geometric setup is much cheaper than the algebraic, especially at the finest levels, a hybrid approach could be cheaper to set up while still allowing accurate restriction of structural features to very coarse grids. This is supported by the performance data in Figure \ref{subfig:2D_Grid_H-BA}, where the hybrid approach uses 3 geometric restriction operators followed by algebraic coarsening to a coarse level with less than 1500 dofs.

In this hybrid approach the setup cost is nearly the same as the GMG approach, but the number of iterations of the linear solver is about the same as the AMG approach, effectively taking the best features from both. This suggests that geometric coarsening can safely be used in any optimization problem up to the point where the minimum dimension of the voids becomes smaller than the dimension of one element on a coarse grid. Beyond that point the effectiveness of GMG preconditioning is dependent on how the features develop. For small problems a pure GMG approach is likely sufficient, but at large scale it should be expected that some elongated features will develop. Using algebraic coarsening may be necessary to prevent merging features over multiple levels of refinement; however, it is difficult to predict at exactly what point the robustness of AMG begins to outweigh the cheap assembly of GMG. This will be covered in more detail in the Discussion section.}

%% file: 3DCantilever.tex
\subsection{3D cantilever beam}
\label{section:3D_Cantilever}
\begin{figure*}
	\centering
	\hspace*{\fill}
	\begin{subfigure}[c]{0.48\textwidth}
		\centering
        \includegraphics[width=\linewidth, clip=True, viewport=200 180 700 540]{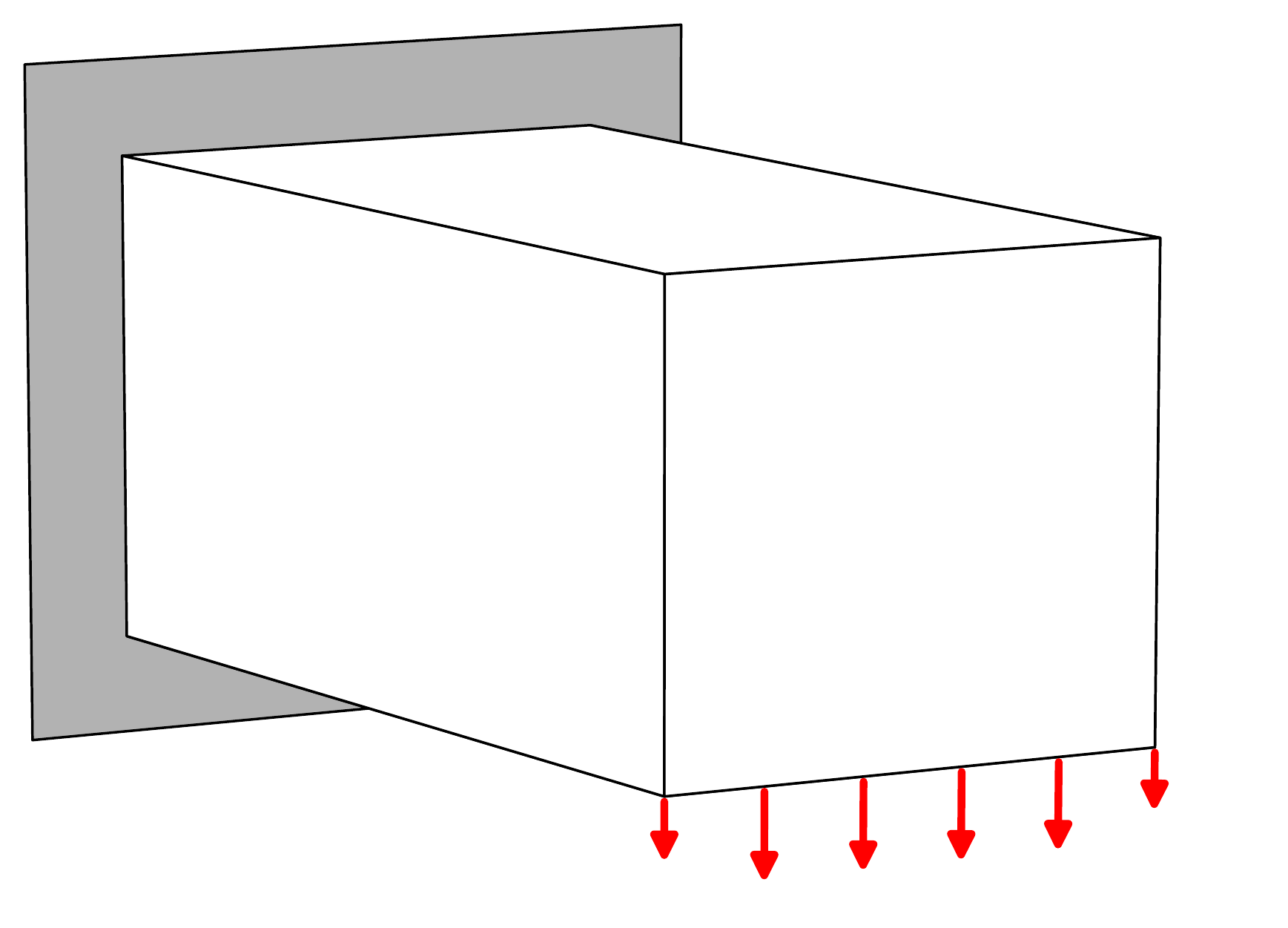}
	\end{subfigure}
    \hfill
	\begin{subfigure}[c]{0.48\textwidth}
        \centering
        \includegraphics[width=0.75\linewidth]{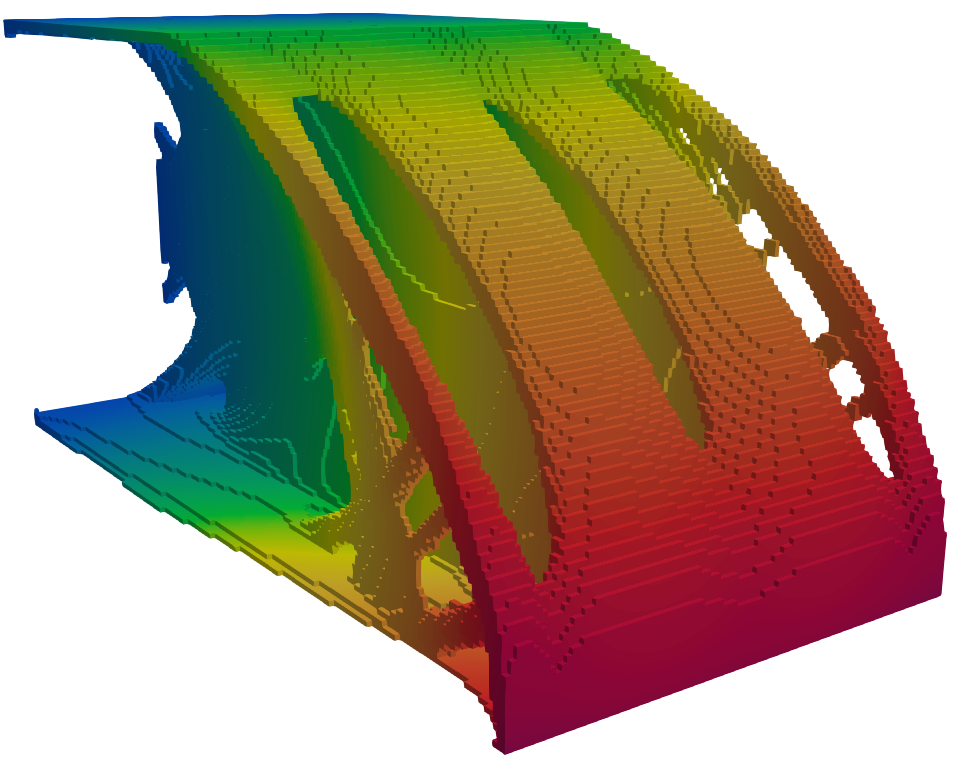}
	\end{subfigure}
	\hspace*{\fill}
    \caption{Design domain and optimized structure for 3D cantilever beam problem. The gray region represents a fixed face of the domain, and a uniform downward load is applied to the bottom of the opposite face. The optimized structure is colored by magnitude of the displacement field.}
    \label{fig:3D_Cantilever}
\end{figure*}

{\add We now present an example describing how the trends observed in 2D translate to 3D optimization problems. We examine the problem of compliance minimization for a cantilever domain in 3D with an aspect ratio of 2:1:1. We use a mesh resolution of 192x96x96 elements across 256 processors and set the volume fraction to 0.12. We perform 20 optimization iterations for each penalty value and the penalty is increased in increments of 0.25. The design domain and result of the highest resolution optimization are shown in Figure \ref{fig:3D_Cantilever}.}
\delete{The next example is compliance minimization for a cantilever domain in 3D with an aspect ratio of 2:1:1. The design domain and result of the highest resolution optimization are shown in Figure \ref{fig:3D_Cantilever}. The problem is run at three different resolutions: 16x8x8 elements, 32x16x16 elements, and 64x32x32 elements. We perform 80 optimization iterations for each penalty value and the penalty is increased in increments of \delete{0.5} {\add 0.25} (the penalty increment is higher than in 2D due to the higher cost of 3D simulations).}

\begin{figure}
	\centering
	\begin{subfigure}[t]{\columnwidth}
    	\centering
        \includegraphics[width=0.85\linewidth]{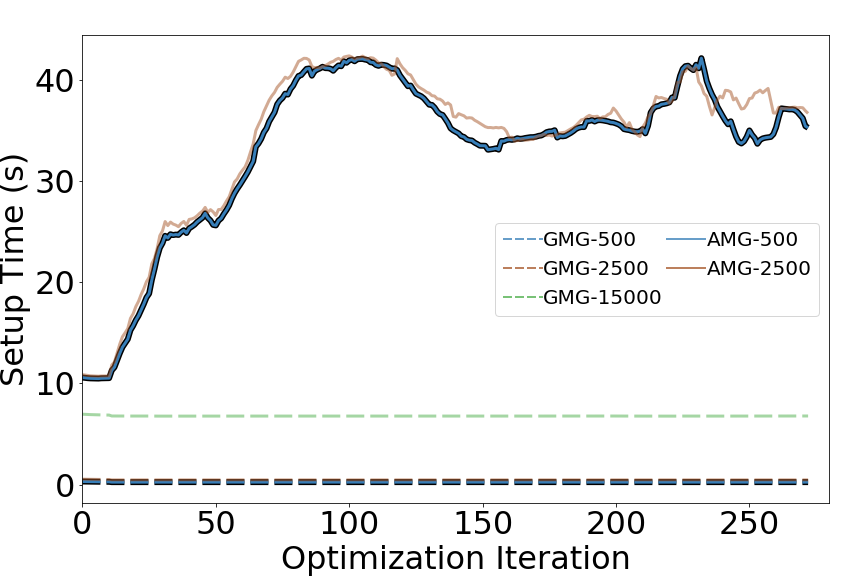}
        \caption{Time to set up the preconditioner.}
    \end{subfigure}

	\begin{subfigure}[t]{\columnwidth}
        \centering
        \includegraphics[width=0.85\linewidth]{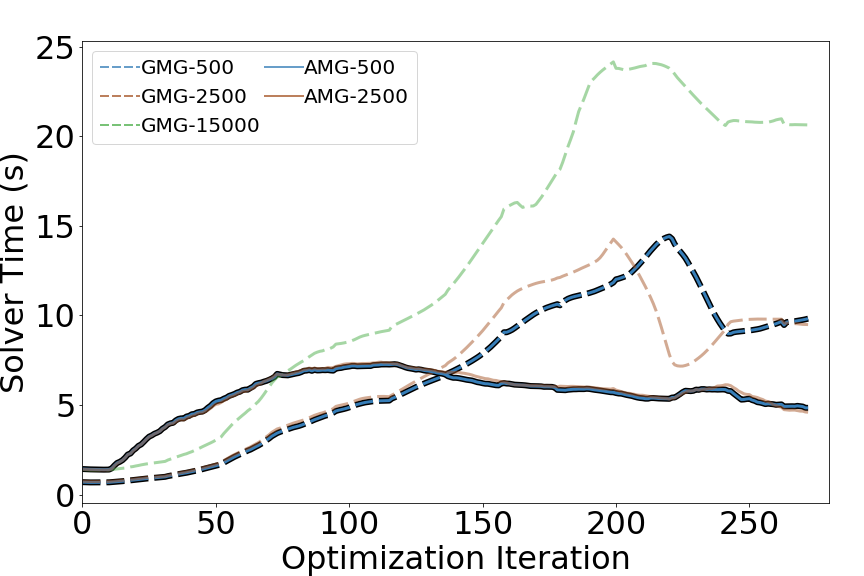}
        \caption{Time to solve the linear system for displacements.}
    \end{subfigure}
    
	\begin{subfigure}[t]{\columnwidth}
        \centering
        \includegraphics[width=0.85\linewidth]{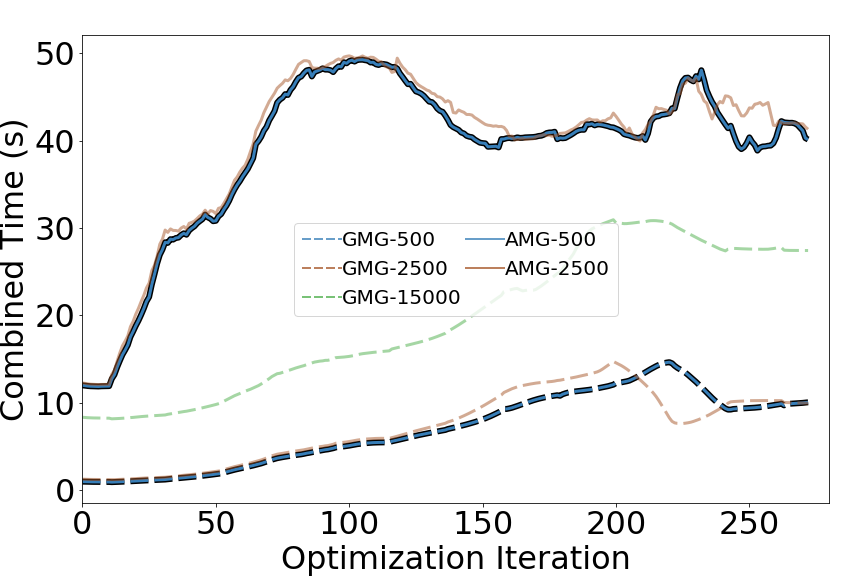}
        \caption{Combined time to set up and solve.}
    \end{subfigure}

    \begin{subfigure}[t]{\columnwidth}
        \centering
        \includegraphics[width=0.85\linewidth]{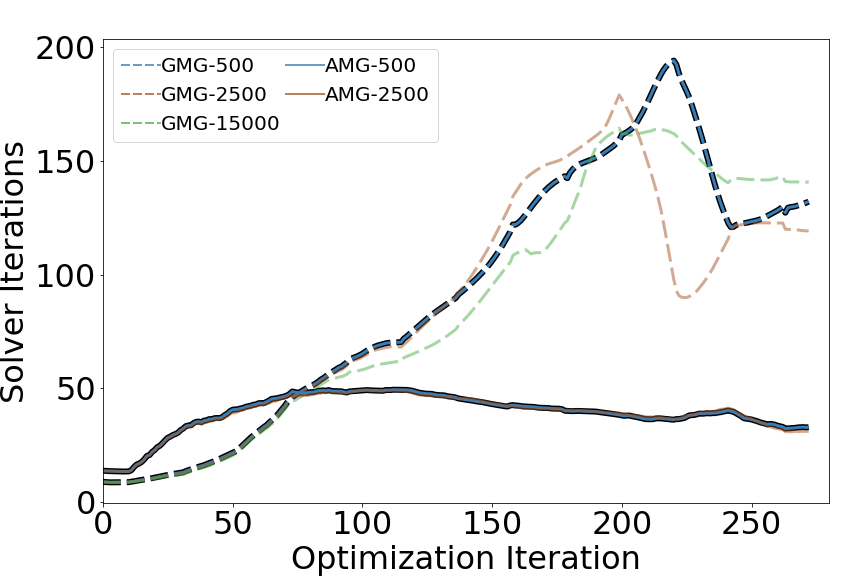}
        \caption{Iterations of the linear solver until convergence.}
    \end{subfigure}
    \caption{Performance of the preconditioners for the 3D cantilever problem with varying coarse grid size. AMG performance is displayed with solid lines and GMG with dashed lines. The best setup for AMG and GMG are emboldened for easier comparison.}
    \label{fig:3D_Cantilever_Coarse_Size}
\end{figure}

{\add As for the 2D problem, we start by using weighted point-block jacobi with a weight of 0.5 as a smoother on all levels, and the coarse grid is solved with an LU decomposition on a single process. Figure \ref{fig:3D_Cantilever_Coarse_Size} demonstrates how the two multigrid versions perform as the coarse grid size is increased. We restrict the coarse grid sizes to fewer than 500, 2500, or 15000 dofs, corresponding to 6, 5, or 4 levels in the geometric hierarchy. The results for the AMG hierarchy with the largest coarse grid size are not shown as the time to set up the preconditioner alone is more than an order of magnitude larger than for the other two coarse grid sizes and the preconditioner as a whole is extremely inefficient compared to the others displayed here.

Similar to before, the setup cost increases with increased size of the coarse grid, and the AMG solver costs much more than GMG. However, contrary to the 2D case, the number of iterations for the GMG preconditioner to converge is largely independent of the coarse grid size. For all three coarse grid sizes the total iterations rarely exceed 200, and as a result increasing the coarse grid size increases the preconditioner cost. The AMG preconditioner similarly requires roughly the same number of iterations regardless of the coarse grid size, though the setup cost is nearly the same for the two coarse grid sizes shown. This is due to the fact that the majority of the setup cost for these cases comes from the matrix triple product instead of the coarse grid factorization. The additional connectivity of nodes in 3D and the larger number of rigid body modes used to construct the prolongation operators means that coarse grid operators are much denser for AMG than GMG, which in turn greatly increases the cost of the preconditioner. In the same vein, setting a large coarse grid size for AMG incurs a huge cost for the coarse grid factorization due to the high density of the operator. These observations mean that for this problem both coarsening strategies are most efficient with very small coarse grid sizes.

Comparing the best setup for either strategy, the setup cost for AMG is much higher than that for GMG, though the number of iterations to convergence is much lower for the AMG preconditioner. While the reduced iterations was enough to offset the setup cost in the 2D case, the greater disparity in setup costs and increased cost of applying the denser AMG operators in 3D results in the GMG preconditioner being much more efficient. In addition, the peak number of iterations in the 3D case is about half that of the 2D case, meaning there is less room for an AMG preconditioner to provide improvement.}

\delete{The overall performance of the various preconditioners on the 3D cantilever problem are shown in Figure \ref{fig:3D_Cantilever_Performance}. The relative performance is similar to the 2D case; however, now the GMG preconditioner outperforms AMG for the smaller two problem sizes. For these problem sizes, the number of iterations to converge is similar between AMG and GMG, both requiring less than 20 iterations at a majority of the optimization steps. For the larger problem size though, as soon as the structure begins to develop, the average number of iterations for the GMG preconditioner jumps to about 75, while the AMG preconditioner only requires about 30 iterations. Again breaking the cost down into the setup and solve phase, the AMG preconditioner takes 2-2.5 times longer to setup than GMG and the block AMG preconditioner takes 2-3 times longer than GMG. In the solve phase the relative performance is drastically different depending on the problem size. For the smaller two problems AMG takes about 40\% longer than GMG, and block AMG only 15\% longer. However, for the larger problem size AMG and block AMG takes 40\% and 50\% less time than GMG, respectively. Combined, the total time to setup and solve the system is 50-80\% longer for either AMG solver in the smaller problems, but about 20\% shorter for the larger problem.}

\begin{figure}
	\centering

	\begin{subfigure}[t]{\columnwidth}
        \centering
        \includegraphics[width=0.85\linewidth]{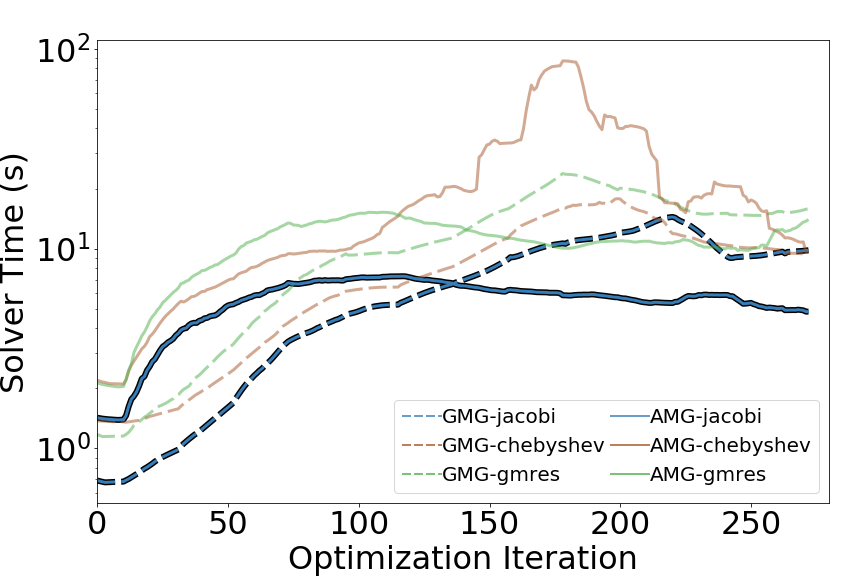}
    \end{subfigure}

    \begin{subfigure}[t]{\columnwidth}
        \centering
        \includegraphics[width=0.85\linewidth]{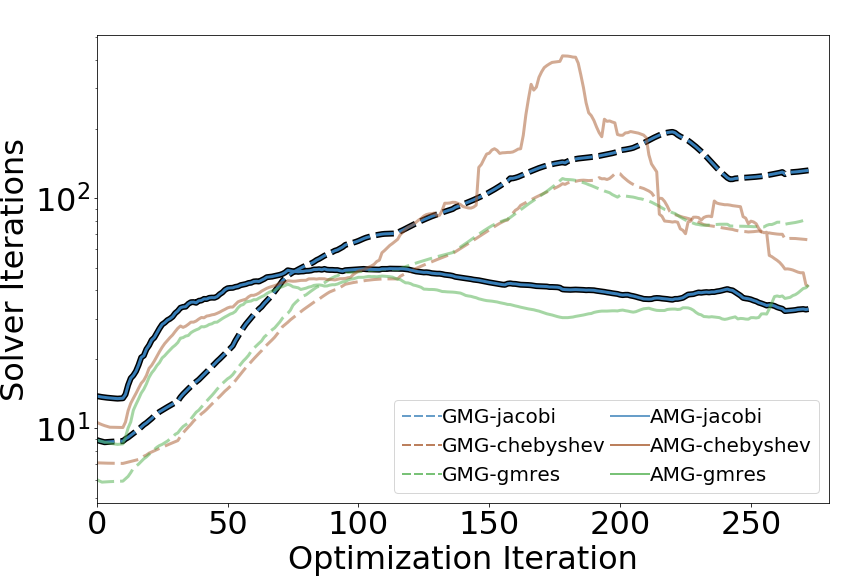}
    \end{subfigure}
    \caption{Performance of the preconditioners for the 3D cantilever problem with different smoothers.}
    \label{fig:3D_Cantilever_Smoothers}
\end{figure}

{\add For completeness, we also show the impact of the different smoothers on the performance of either multigrid setup in Figure \ref{fig:3D_Cantilever_Smoothers}. All of the same patterns observed in the 2D case appear again in 3D. Notably, the SOR-GMRES smoother requires the fewest iterations, the weighted Jacobi smoother requires the least time, and the SOR-Chebyshev smoother experiences intermittent stagnation. We again conclude that the weighted Jacobi smoother can be regarded as the most reliable in this case.}

\begin{figure}
	\centering
	\begin{subfigure}[t]{\columnwidth}
    	\centering
        \includegraphics[width=0.85\linewidth]{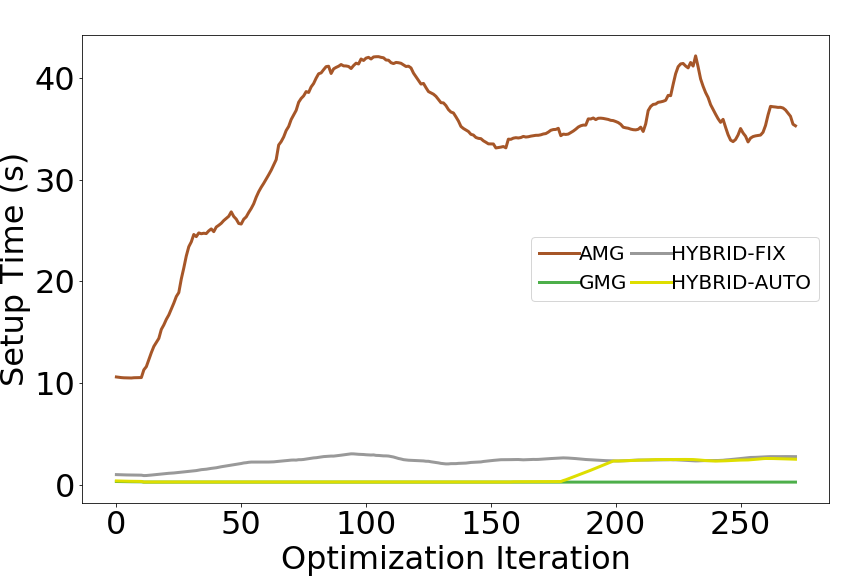}
        \caption{Time to set up the preconditioner.}
    \end{subfigure}

	\begin{subfigure}[t]{\columnwidth}
        \centering
        \includegraphics[width=0.85\linewidth]{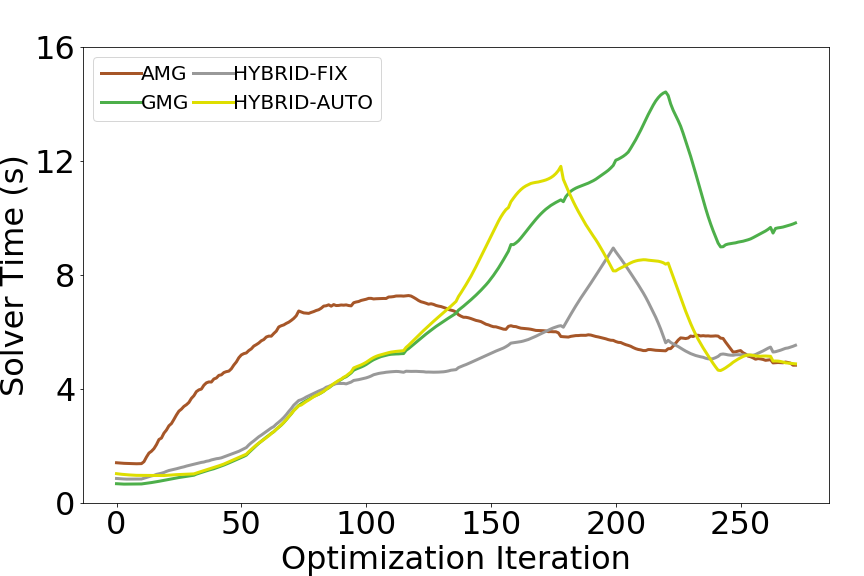}
        \caption{Time to solve the linear system for displacements.}
    \end{subfigure}
    
	\begin{subfigure}[t]{\columnwidth}
        \centering
        \includegraphics[width=0.85\linewidth]{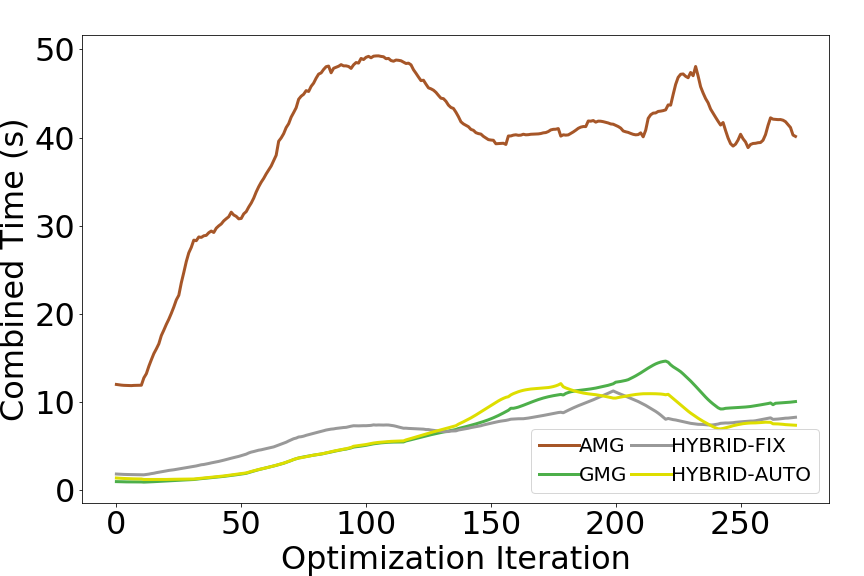}
        \caption{Combined time to set up and solve.}
    \end{subfigure}

    \begin{subfigure}[t]{\columnwidth}
        \centering
        \includegraphics[width=0.85\linewidth]{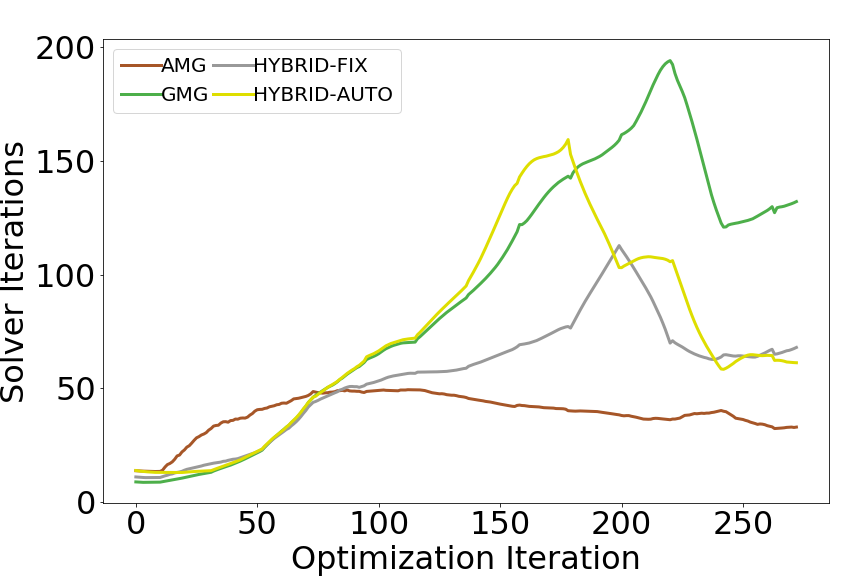}
        \caption{Iterations of the linear solver until convergence.}
    \end{subfigure}
    \caption{Comparison of AMG, GMG, and Hybrid performance for the 3D cantilever problem.}
    \label{fig:3D_Cantilever_Hybrid_MG}
\end{figure}

{\add Recognizing that the biggest disadvantage of the AMG preconditioner in 3D is the high setup cost, we also compare the performance of \deleted{a}{\addd two different} hybrid GMG-AMG preconditioner{\addd s} to the optimal GMG and AMG preconditioners in Figure \ref{fig:3D_Cantilever_Hybrid_MG}. {\addd In the first approach we} \deleted{We} restrict the coarse grid of the hybrid preconditioner to have less than 1000 dofs and start the hierarchy with 2 geometric restriction operators. {\addd In the second approach we restrict the coarse grid to the same size, but we start with 6 GMG levels instead of 2. We then reduce the number of GMG levels by 1 every time the linear solver requires more than 200 iterations to converge, until a minimum of 2 GMG levels remain.}  As the figure shows, \deleted{this} {\addd both approaches} drastically reduce\deleted{s} the setup cost of the preconditioner when compared to a pure AMG setup, while keeping the number of iterations lower than the pure GMG preconditioner. Combining these properties allows the hybrid preconditioner{\addd s} to keep the total time to setup and solve substantially lower than the pure AMG preconditioner and even slightly lower than the pure GMG preconditioner on average.}

{\addd The second hybrid MG approach provides additional improvement in performance by only using as many AMG levels as necessary. When geometric restriction is sufficient in keeping iteration counts low, it is favored for it's much cheaper setup and application costs (due to increased operator sparsity). As performance of the geometric approach degrades, the preconditioner progressively transitions to the more expensive algebraic approach to maintain a low iteration count. This results in a 3\% reduction of total time to setup the preconditioner and solve the linear system across the entire optimization process for this particular example. Noting that the automated hybrid approach experiences significantly higher iteration counts from optimization steps 100 through 200, it appears that more aggressively transitioning from GMG to AMG could yield greater performance improvement, though we do not go so far as to describe an optimal transitioning strategy here.}

%% file: Discussion.tex
\section{Discussion}
\label{section:Discussion}
We have analyzed the relative performance of geometric (GMG) and algebraic (AMG) preconditioners in the context of topology optimization. For topology optimization at large scales it is necessary to use iterative solvers, which rely on effective preconditioners for their performance. AMG and GMG preconditioners both use the same basic procedure for solving or preconditioning a linear system; however, the methods differ in how the preconditioners are constructed. In GMG, the coarser grids are constructed directly from the mesh that discretizes the problem domain, ignoring the evolution of structural features throughout the design optimization. In contrast, AMG methods perform grid coarsening based on the stiffness matrix alone, without any direct knowledge of the underlying discretization of the partial differential equation (PDE).

In all of our examples we have seen that the AMG preconditioner is more expensive to construct due to the nature of the interpolation/restriction operators. \delete{However,} {\add For 2D problems} the setup cost is generally offset by the fact that the AMG preconditioner more actively adapts to changing structural topology. For some simple cases seen often in the literature, the performance of GMG preconditioners is similar to \delete{or slightly better than} {\add that of} AMG preconditioners. However, we have also demonstrated \delete{several cases} {\add criteria} where the AMG preconditioners are much more robust due to the extra work in their assembly. The increased robustness comes from the inherent capacity for AMG to identify where structural features exist while constructing the hierarchy. In addition, AMG preconditioners are much \delete{better suited}{\add simpler to employ} for optimization problems on irregular meshes, for example when polygonal meshes are used \citep{Talischi2012}.

{\add In 3D the AMG preconditioners are much less competitive, mainly due to the high cost of assembly. The increased nodal connectivity of 3D meshes, as well as the larger number of rigid body modes used for smoothed aggregation interpolation, lead to much denser operators in the AMG hierarchy. Some strategies exist to alleviate these issues, such as more aggressive graph coarsening, but they come at the cost of less efficient preconditioning. In the experience of the authors they aren't able to reduce the overall cost of the AMG preconditioner enough to make it competitive. We conclude that it is very difficult to make AMG more efficient than GMG for 3D optimization problems, though it is unclear if this holds for scalar-valued problems, such as thermal conduction, where a smaller near-nullspace (1 dimension instead of 6) leads to smaller and cheaper coarse grid operators.

As an alternative to using AMG, particularly in 3D where the cost is high, we also examine the merit of using a hybrid GMG-AMG preconditioner as described in \citep{Aage2017Giga-voxelDesign}. We show that using geometric restriction for the finest levels of the hierarchy drastically reduces the setup cost compared to pure AMG, while algebraic restriction at the coarser levels helps to keep the linear solver iterations consistently lower than for a pure GMG preconditioner. The level of coarsening where the restriction should change from geometric to algebraic is an open question, though in the author's experience using an approximately equal amount of geometric and algebraic restriction consistently produces good results. \deleted{A possible adaptive approach would be to switch the coarsest restriction methods to algebraic as the linear solver iterations increase. In all the examples presented here, iterations of the solver when using AMG almost never exceed 100. Therefore, adaptively switching coarse levels from geometric to algebraic restriction when the number of iterations exceed 100-200 can also be expected to produce good results.} {\addd We also describe an approach to adaptively change between GMG and AMG to further improve performance as the structure evolves. Noting that AMG is more expensive per iteration, but more effective at reducing the total number of iterations, we progressively reduce the number of geometric levels in the MG hierarchy as the solver iteration count rises. This offers a further advantage of cheap GMG coarsening when iteration counts are low and more effective (though more expensive) AMG coarsening as performance deteriorates. An application of this adaptive approach to the 3D cantilever problem at a larger scale is described in the appendix.}

Another important conclusion is that GMG smoothers work very well as long as they don't prematurely lump structural features together. For topologies where all features have similar aspect ratios (this could also be framed as having voids with small aspect ratios), this is rarely a problem and GMG is likely a better choice of preconditioner. If these conditions are not satisfied, AMG coarsening is necessary to keep iterations of the iterative solver low. However, AMG coarsening is much more expensive, particularly in 3D. To that end, we recommend that hybrid GMG-AMG preconditioners be used, particularly in 3D, to take advantage of cheap coarsening at the finer levels of the hierarchy while using more robust algebraic coarsening at the coarser levels. This keeps both the setup cost of the preconditioner and iterations of the solver low and will likely be necessary to run problems at extreme scales. A hybrid approach can also be beneficial in 2D, though the cheaper cost of creating and using coarse scale AMG operators means that the benefit is small. Given that a black-box AMG preconditioner is also simpler to use, we recommend that approach for 2D problems.}

The relative merit of GMG or AMG preconditioners is also highly dependent on the type of optimization being performed. In the simple case of compliance minimization, where only a single solution to the linear system is needed, the GMG preconditioner is \delete{often just as effective as the AMG preconditioner}{\add less robust than the AMG preconditioner in reducing the number of iterations, but very competitive in terms of runtime}. In cases where additional solutions to the linear system are needed, for example to solve adjoint problems or evaluating structural stability using the generalized eigenvalue problem, AMG methods are \delete{much}{\add relatively} more effective. In these cases the better performance of the iterative solver with an AMG preconditioner is more important than the cheaper setup cost of the GMG preconditioner.

{\add We have also compared the relative performance of 3 popular smoothing choices already used in the topology optimization literature. It is our opinion that weighted Jacobi smoothing is the most robust, however the performance benefit is slight when compared to other factors, such as coarse grid size. In addition to the modest performance benefit, it is also a simple preconditioner to implement and tune. While SOR-preconditioned Chebyshev iterations are well-regarded for smoothing, we have demonstrated that they are less robust than weighted Jacobi for the examples presented here.}

We conclude that \delete{for simpler optimization formulations such as compliance minimization (only one solution to the linear system is needed) on uniform grids,}the relative simplicity of GMG makes it a very appealing preconditioner{\add, especially in 3D}. However, for problems where multiple solutions to the same system are needed or the topology exhibits certain characteristics, such as tightly-packed, highly flexible features, AMG offers a substantial performance improvement. In addition, AMG readily extends to problems on irregular domains or non-uniform meshes, and may be combined with other cost saving measures, such as design space optimization, to further improve performance. {\add GMG is not inherently incompatible with any of these scenarios, but is simply much more complicated to implement, whereas AMG requires no extra efforts from the user.} Algebraic multigrid has also experienced sufficient development that near black-box functions are available in most scientific computing environments \citep{petsc-user-ref,OlSc2018,Falgout2002Hypre:Preconditioners}. {\add In that sense, while the actual work of constructing an AMG precondtioner is more complicated than GMG, it is likely easier than GMG to implement through one of these black box solvers\deleted{,}\addd{.}}


%% file: Appendix.tex
\begin{figure*}
	\centering
	\hspace*{\fill}
	\begin{subfigure}[c]{0.45\textwidth}
        \centering
        \includegraphics[width=0.8\linewidth]{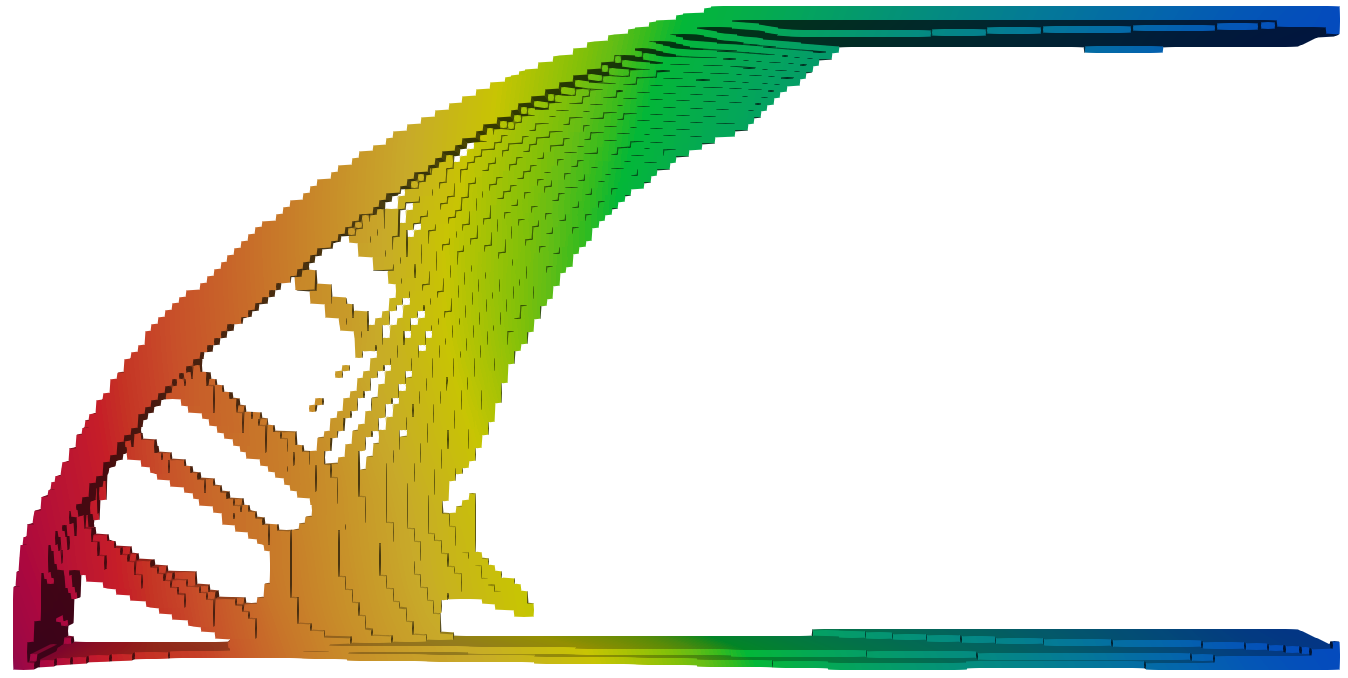}
	\end{subfigure}
    \hfill
	\begin{subfigure}[c]{0.45\textwidth}
        \centering
        \includegraphics[width=0.8\linewidth]{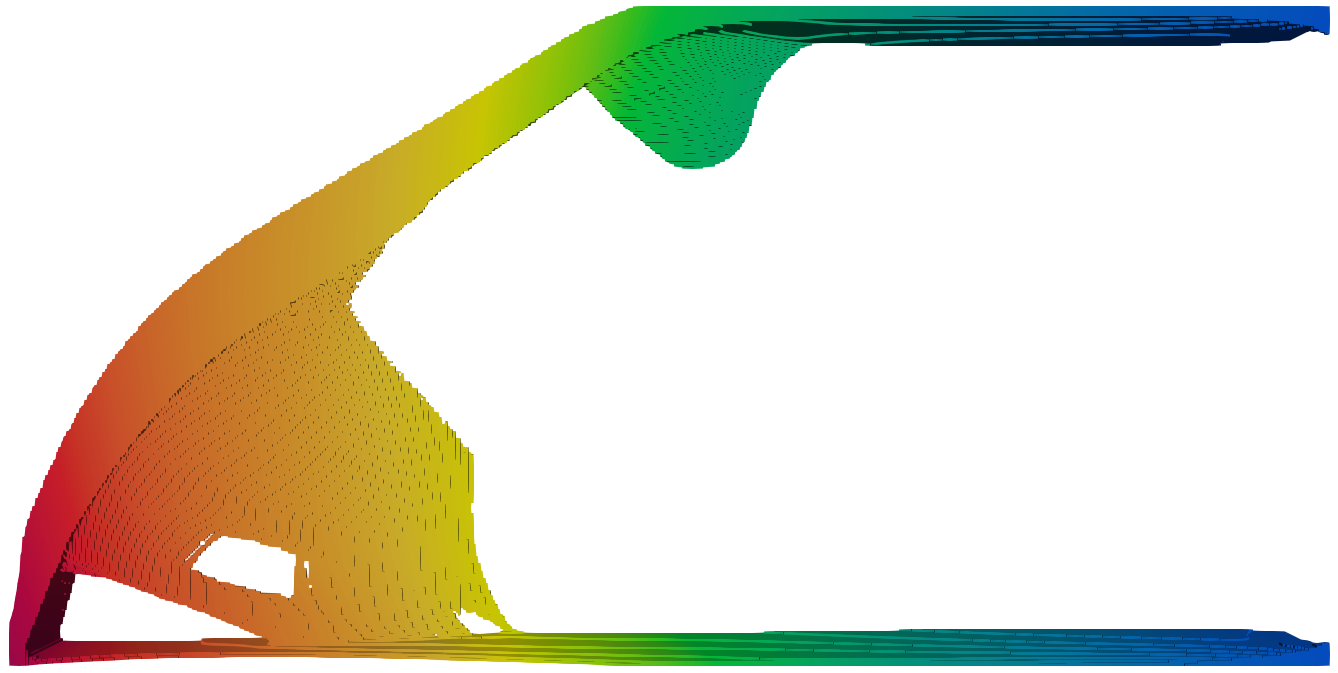}
	\end{subfigure}
	\hspace*{\fill}
    \caption{Comparison of structure along lateral edge of domain at coarse (left) and fine (right) resolution.}
    \label{fig:3D_Cantilever_Side}
\end{figure*}
\begin{figure*}
	\centering
	\hspace*{\fill}
	\begin{subfigure}[c]{0.45\textwidth}
        \centering
        \includegraphics[width=0.4\linewidth]{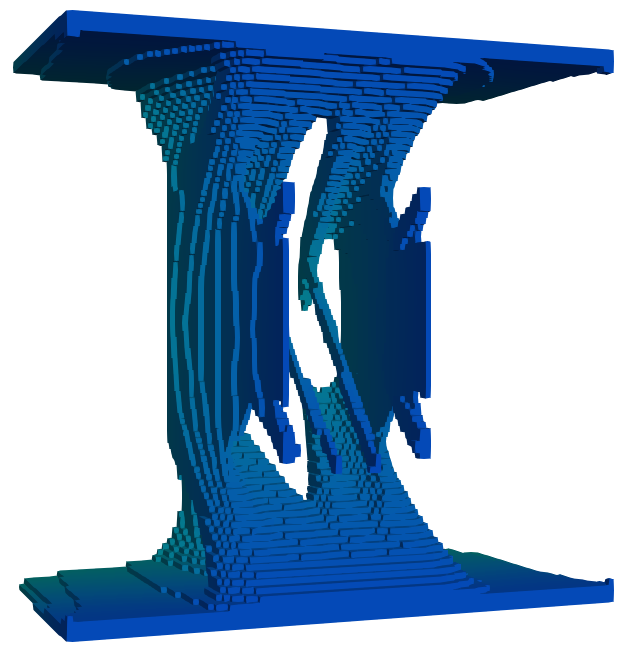}
	\end{subfigure}
    \hfill
	\begin{subfigure}[c]{0.45\textwidth}
        \centering
        \includegraphics[width=0.4\linewidth]{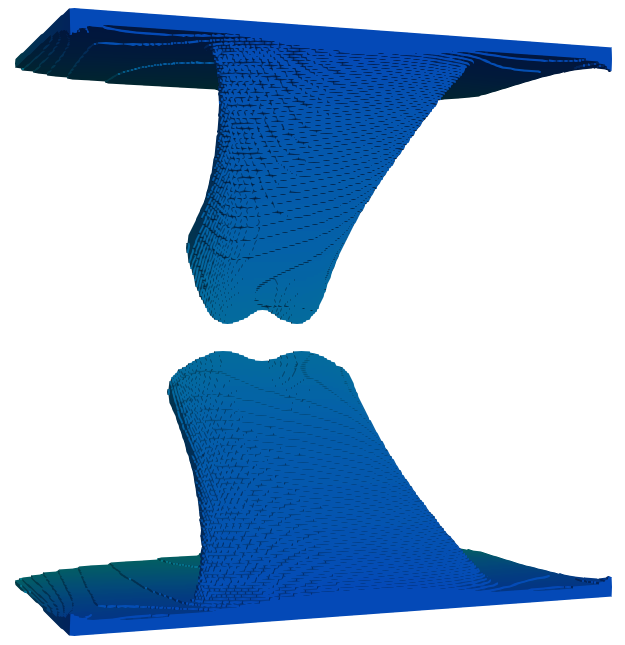}
	\end{subfigure}
	\hspace*{\fill}
    \caption{Comparison of structure along rear edge of domain at coarse (left) and fine (right) resolution.}
    \label{fig:3D_Cantilever_Back}
\end{figure*}

{\addd
\section{Large Scale 3D cantilever beam}
\label{section:3D_Cantilever_Large}
Here we describe the application of our proposed hybrid GMG-AMG preconditioner to the 3D cantilever problem at a larger scale (512x256x256 element mesh, approximately 100e6 dofs). This time we only perform 16 optimization iterations for each penalty increment to improve runtime, but the penalty is still increased in increments of 0.25 to a final value of 4. The filter radius is again set to 1.5 times the element dimension, meaning that the physical dimension of the filter radius is reduced by a factor of 2.67. The final result is shown in Figure \ref{fig:3D_Cantilever_Large}.

\begin{figure*}
	\centering
	\hspace*{\fill}
	\begin{subfigure}[c]{0.45\textwidth}
        \centering
        \includegraphics[width=0.6\linewidth]{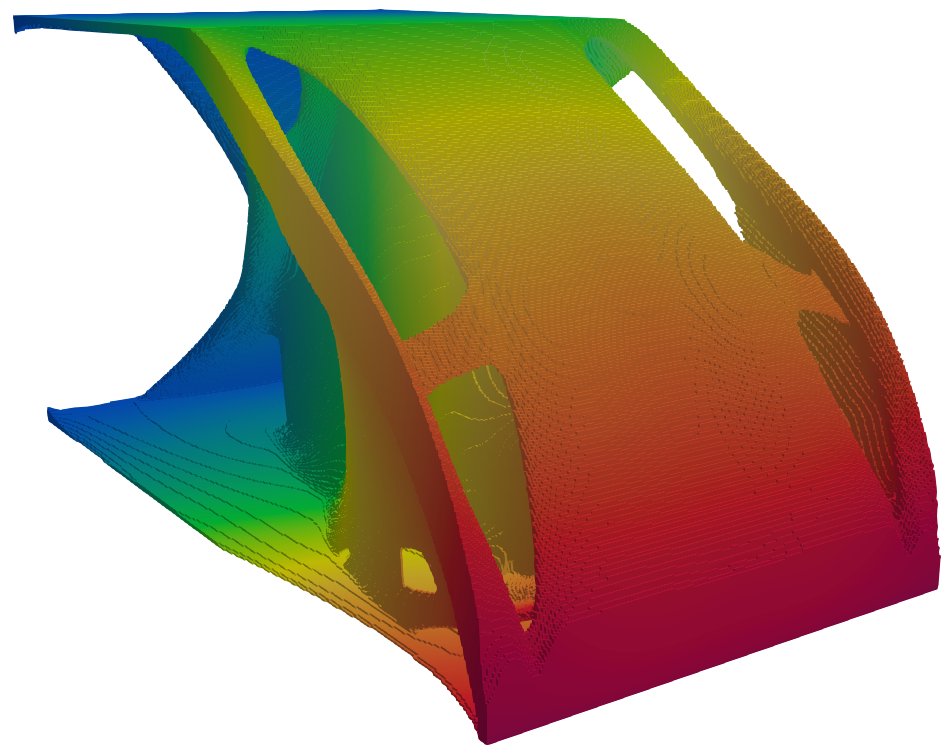}
	\end{subfigure}
    \hfill
	\begin{subfigure}[c]{0.45\textwidth}
        \centering
        \includegraphics[width=0.6\linewidth]{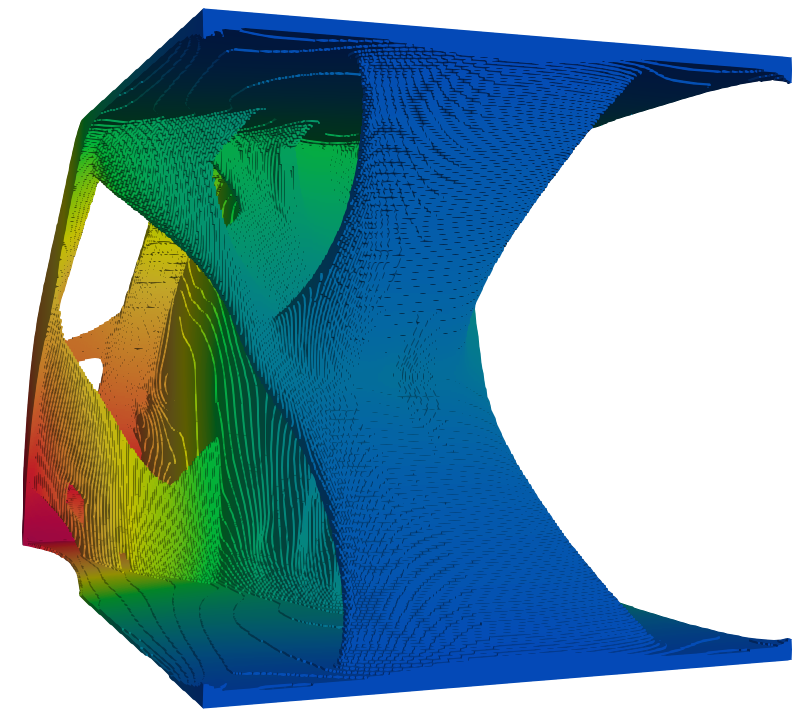}
	\end{subfigure}
	\hspace*{\fill}
    \caption{Optimized structure for 3D cantilever beam problem at an increased resolution.}
    \label{fig:3D_Cantilever_Large}
\end{figure*}

In Figure \ref{fig:3D_Cantilever_Large_Performance} we show the performance of using the hybrid preconditioner for this problem. We restrict the coarse grid size to be less than 1000 dofs, corresponding to 7 levels in the geometric hierarchy, and we apply weighted point-block jacobi as the smoother. It is important to note that in this case, the number of iterations to convergence remains substantially lower than in the smaller problem, never reaching even 100 iterations in a single optimization step. As a result, the hybrid scheme never converts levels of the multigrid preconditioner from geometric to algebraic coarsening. This means that 7 geometric levels are used through the entire simulation, with one small algebraic level at the end for ease of implementation.

\begin{figure*}
	\centering
	\begin{subfigure}[t]{\columnwidth}
    	\centering
        \includegraphics[width=0.85\linewidth]{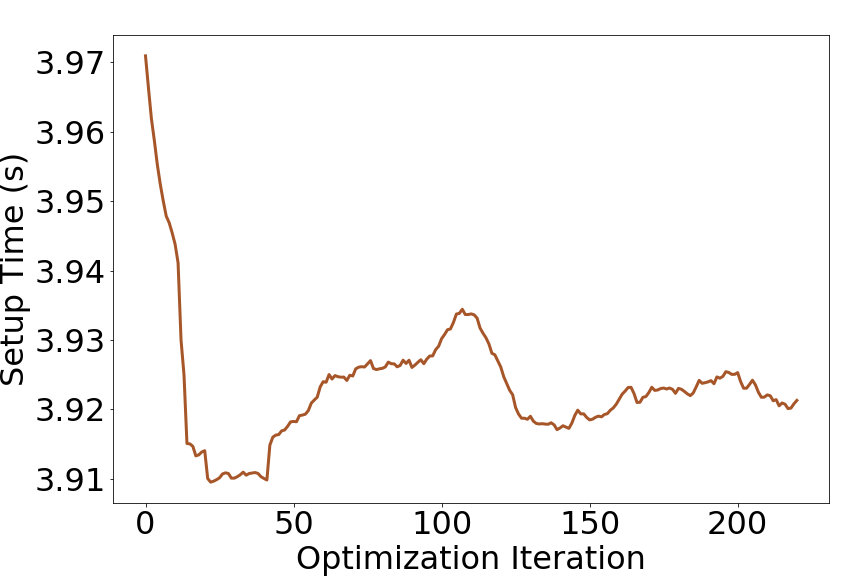}
        \caption{Time to set up the preconditioner.}
    \end{subfigure}
    \hfill
	\begin{subfigure}[t]{\columnwidth}
        \centering
        \includegraphics[width=0.85\linewidth]{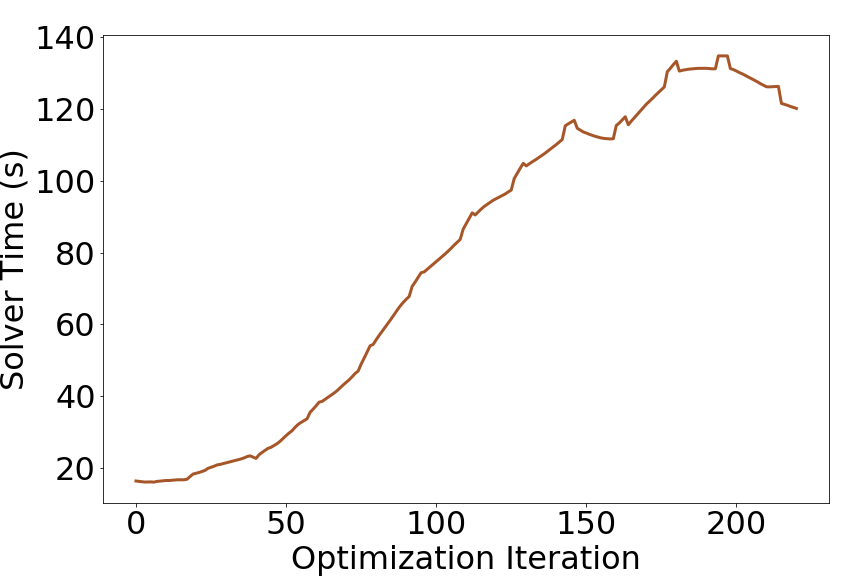}
        \caption{Time to solve the linear system for displacements.}
    \end{subfigure}
    
	\begin{subfigure}[t]{\columnwidth}
        \centering
        \includegraphics[width=0.85\linewidth]{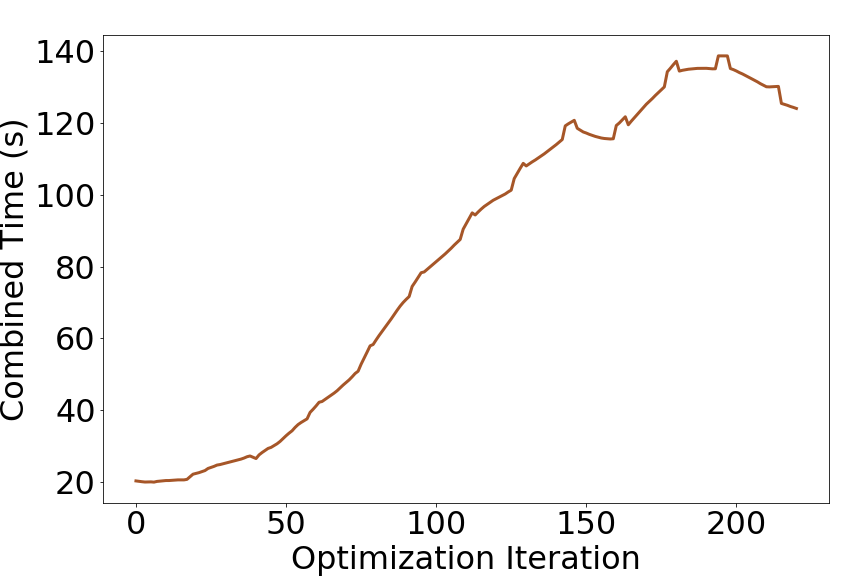}
        \caption{Combined time to set up and solve.}
    \end{subfigure}
    \hfill
    \begin{subfigure}[t]{\columnwidth}
        \centering
        \includegraphics[width=0.85\linewidth]{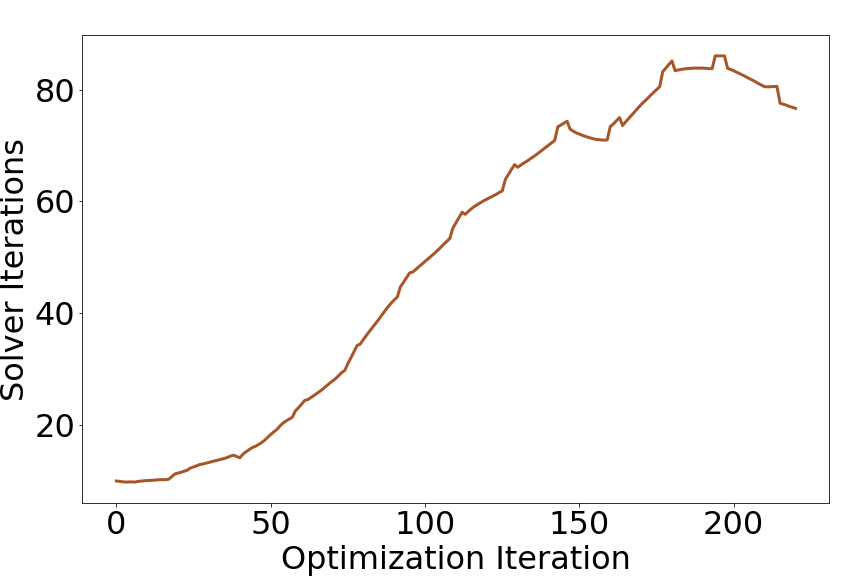}
        \caption{Iterations of the linear solver until convergence.}
    \end{subfigure}
    \caption{Performance of the hybrid preconditioner on the large 3D cantilever problem.}
    \label{fig:3D_Cantilever_Large_Performance}
\end{figure*}

It could be argued that this is an indication that AMG is not needed in larger scale problems; however, an important observation arises upon closer inspection of the optimized structures. Figure \ref{fig:3D_Cantilever_Side} shows the structure that develops along one edge of the domain for both the smaller and larger 3D problem and Figure \ref{fig:3D_Cantilever_Back} shows the difference in the structures along the rear support. Note that in the small case, many more small structural elements develop in close proximity to each other in both of these regions. It is exactly these types of elements that cause difficulty for the GMG preconditioner, but they do not develop in the larger problem. Instead, these numerous structural elements are replaced with a single smooth feature that can only be captured on the higher-fidelity design space. As it is impossible to predict a priori what type of structure will develop from a given loading condition and mesh resolution this motivates the use of our hybrid preconditioner, which only applies algebraic coarsening as necessary. In the smaller problem the algebraic coarsening improves performance when geometric coarsening struggles, but in the larger problem the cheaper geometric strategy is used. In either case there is no need for the user to specify a fixed number of algebraic or geometric levels in the multigrid hierarchy.}

